\documentclass{article}
\usepackage{amsmath}
\usepackage{amsfonts}
\usepackage{amssymb}	
\usepackage{color}
\newtheorem{theo}{Theorem}[section]
\newtheorem{lem}[theo]{Lemma}
\newtheorem{cor}[theo]{Corollary}
\newtheorem{rem}[theo]{Remark}
\newtheorem{prop}[theo]{Proposition}
\newtheorem{defi}{Definition}[section]

\newcommand{\mysection}[1]{\section{#1} \setcounter{equation}{0}}
\newcommand{\proof}{{\sc Proof.} \quad}

\newcommand{\R}{\mathbb{R}}
\newcommand{\N}{\mathbb{N}}
\newcommand{\be}{\begin{equation} \label}
\newcommand{\ee}{\end{equation}}
\newcommand{\bes}{\begin{equation} \begin{array}{c} \label}
\newcommand{\ees}{\end{array} \end{equation}}
\newcommand{\bea}{\begin{eqnarray}\label}
\newcommand{\eea}{\end{eqnarray}}
\newcommand{\beas}{\begin{eqnarray} \begin{array}{rcl} \label}
\newcommand{\eeas}{\end{array} \end{eqnarray}}
\newcommand{\bas}{\begin{eqnarray*}}\newcommand{\eas}{\end{eqnarray*}}
\newcommand{\bass}{\begin{eqnarray*} \begin{array}{rcl}}
\newcommand{\eass}{\end{array} \end{eqnarray*}}
\newcommand{\basss}{\begin{eqnarray*} \begin{array}{c}}
\newcommand{\easss}{\end{array} \end{eqnarray*}}
\newcommand{\qed}{{}\hfill $\square$ \\}
\newcommand{\bit}{\begin{itemize}}
\newcommand{\eit}{\end{itemize}}
\newcommand{\nn}{\nonumber}
\newcommand{\eps}{\varepsilon}
\newcommand{\abs}{\\[3mm]}

\newcommand{\pO}{\partial\Omega}

\newcommand{\io}{\int_\Omega}

\newcommand{\F}{{\cal F}}
\newcommand{\D}{{\cal D}}
\newcommand{\B}{{\cal B}}
\newcommand{\set}{{\cal S}(m,M,B,\kappa)}
\newcommand{\tm}{T_{max}(u_0,v_0)}

\begin{document}
\title{New critical exponents in a fully parabolic quasilinear Keller-Segel and applications to volume filling models}
\author{
Tomasz Cie\'{s}lak \\
{\small Instytut Matematyczny PAN, \'Sniadeckich 8, 00-956 Warsaw, Poland}\\
{\small E-Mail: T.Cieslak@impan.pl}
\and
Christian  Stinner\\
{\small Felix-Klein-Zentrum f\"{u}r Mathematik, TU Kaiserslautern, Paul-Ehrlich-Str. 31,}\\ {\small 67663 Kaiserslautern, Germany.
E-Mail: stinner@mathematik.uni-kl.de}
}
\date{}
\maketitle
\begin{abstract}
  \noindent
 We carry on our studies related to the fully parabolic quasilinear Keller-Segel system started in \cite{Cie_Sti_JDE}
 and continued in \cite{Cie_Sti_AAM}. In the above mentioned papers we proved finite-time blowup of radially symmetric
 solutions to the quasilinear Keller-Segel system if the nonlinear chemosensitivity is strong enough and an adequate
 relation between nonlinear diffusion and chemosensitivity holds. On the other hand we proved that once
 chemosensitivity is weak enough solutions exist globally in time. The present paper is devoted to looking for critical
 exponents distinguishing between those two behaviors. Moreover, we apply our results to the so-called volume filling
 models with a power-type probability jump function. The most important consequence of our investigations of the latter
 is a critical mass phenomenon found in dimension $2$. Namely we find a value $m_*$ such that when the solution to the two-dimensional volume filling Keller-Segel system starts with mass smaller than $m_*$, then it is bounded, while for initial data with mass exceeding $m_*$ solutions are unbounded, though being defined for any time $t>0$.


\end{abstract}
\mysection{Introduction}\label{section1}
This work continues our investigations related to nonnegative solution couples $(u,v)$ of the parabolic-parabolic
Keller-Segel system
\be{0}
	\left\{ \begin{array}{ll}
	u_t= \nabla \cdot (\phi(u) \nabla u) - \nabla \cdot ( \psi(u) \nabla v ), & \; x\in\Omega, \ t>0, \\[2mm]
	v_t=\Delta v-v+u, & \; x\in\Omega, \ t>0, \\[2mm]
	\frac{\partial u}{\partial\nu}=\frac{\partial v}{\partial\nu}=0, & \; x\in\partial\Omega, \ t>0, \\[2mm]
	u(x,0)=u_0(x), \quad v(x,0)=v_0(x), & \; x\in\Omega,
	\end{array} \right.
\ee
in $\Omega \subset \R^n$, where $n \ge 2$, and the
initial data are supposed to satisfy $u_0 \in C^0(\bar\Omega)$ and $v_0\in W^{1,\infty}(\Omega)$ such that
$u_0 \geq 0$ and $v_0 \geq 0$ in $\bar{\Omega}$.

Moreover, we assume that $\phi, \psi \in C^2([0,\infty))$ and that there is $\beta \in C^2 ([0,\infty))$
such that
\begin{equation}\label{0.1}
 \phi(s) >0, \qquad \psi(s) = s \beta (s), \quad\mbox{and }\quad \beta (s) >0 \quad\mbox{for } s \in [0,\infty)
\end{equation}
are satisfied.

In previous papers \cite{Cie_Sti_JDE} and \cite{Cie_Sti_AAM} we have shown the finite-time blowup of solutions to \eqref{0} in dimensions $n\geq 3$ and $n=2$ respectively provided $\psi$ is superlinear and a proper relation between
$\psi$ and $\phi$ holds. It was known already that if the above mentioned relation between $\psi$ and $\phi$ is not satisfied then solutions exist globally in time (see for instance \cite{taowin_jde}). However, the question arose whether the superlinear condition on $\psi$ is necessary. We managed to show that actually some restriction on the growth of $\psi$ has to be imposed. Otherwise, even in the case of parameters yielding finite-time blowup for
$\psi$ not decaying, when $\psi$ decays sufficiently fast, solutions to \eqref{0} must exist globally in time. However, we indicated the examples of unbounded radially symmetric solutions in that case. The aim of the present paper is to find critical exponents on the growth of $\psi$ distinguishing between possibility of finite-time blowup and the lack of it when $\psi$ and $\phi$ satisfy the supercritical relation.

Before moving to the presentation of the precise results, let us mention that the finite-time blowup results for the fully parabolic Keller-Segel system were unavailable for years, the only existing result in the literature before 2010 being \cite{HV:1}, where a difficult construction of a nongeneric special example of radially symmetric blowing up solution of a semilinear two-dimensional Keller-Segel system was given. Next results appeared only recently, first
finite-time blowup for large mass solutions of a one-dimensional supercritical quasilinear Keller-Segel system, see \cite{clKS}, next the breakthrough due to M.Winkler in \cite{win_bu} and our considerations being an extension of the method of Winkler to the quasilinear case in \cite{Cie_Sti_JDE, Cie_Sti_AAM}. For a more detailed discussion on this issue and the discussion on the known results for parabolic-elliptic Keller-Segel we refer the interested reader to the Introduction in \cite{Cie_Sti_JDE}.

Suppose that there exist positive constants $s_0$, $a$, and $b$ such that the functions
\begin{equation}\label{GH}
 G(s) := \int\limits_{s_0}^s \int\limits_{s_0}^\sigma \frac{\phi(\tau)}{\psi(\tau)} \; d\tau \, d \sigma, \quad s>0,
 \qquad\mbox{and}\qquad
 H(s) := \int\limits_0^s \frac{\sigma \phi(\sigma)}{\psi(\sigma)} \; d\sigma, \quad s \ge 0,
\end{equation}
fulfill
\begin{equation}\label{G1}
  G(s) \le a \, s^{2-\alpha}, \; s \ge s_0,  \qquad\mbox{with some } \alpha > \frac{2}{n},
\end{equation}
as well as
\begin{equation}\label{H1}
  H(s) \le \gamma \cdot G(s) + b(s+1), \; s > 0, \qquad\mbox{with some } \gamma \in \left( 0, \frac{n-2}{n}
  \right).
\end{equation}
We remark that $H$ in \eqref{GH} is well-defined due to the positivity of $\beta$ in $[0,\infty)$.

It is well-known that the function
\begin{equation}\label{F}
  \F(u,v):=\frac{1}{2} \io |\nabla v|^2 + \frac{1}{2} \io v^2 - \io uv + \io G(u)
\end{equation}
is a Liapunov functional for \eqref{0} with dissipation rate
\begin{equation}\label{D}
  \D(u,v):=\io v_t^2 + \io \psi(u) \cdot \Big| \frac{\phi (u)}{\psi(u)} \nabla u - \nabla v \Big|^2.
\end{equation}
More precisely, any classical solution to \eqref{0} satisfies
\begin{equation}\label{liapunov}
	\frac{d}{dt} \F(u(\cdot,t),v(\cdot,t)) = - \D (u(\cdot,t),v(\cdot,t))
	\qquad \mbox{for all } t \in (0,\tm),
\end{equation}
where $\tm \in (0,\infty]$ denotes the maximal existence time of $(u,v)$ (see \cite[Lemma~2.1]{win_mmas}). \abs

Assume further that
\begin{equation}\label{psi}
  \frac{s^2}{\psi(s)} \le L \, ( G(s) +s+1), \quad s > 0,
\end{equation}
is satisfied with some $L>0$. This condition enables us to extend our previous finite-time blowup results to a wider class of functions $\psi$.  Namely, we have the following result.
\begin{theo}\label{theo1}
  Suppose that $\Omega=B_R \subset \R^n$ with some $n\ge 3$ and $R>0$, assume that \eqref{G1}, \eqref{H1}, and \eqref{psi}
  are satisfied, and let $m>0$ and $A>0$ be given.
  Then there exist positive constants $T(m,A)$ and $K(m)$ such that for any
  \begin{eqnarray}\label{t1.1}
	(u_0,v_0) \in \B (m,A) &:=& \bigg\{
	(u_0,v_0) \in C^0(\bar\Omega) \times W^{1,\infty}(\Omega) \ \bigg| \
	\mbox{$u_0$ and $v_0$ are radially symmetric} \nn \\
    & & \hspace*{5mm} \mbox{and positive in $\bar\Omega$,
	$\io u_0=m$, $\|v_0\|_{W^{1,2}(\Omega)} \le A$,}\nn \\
    & & \hspace*{5mm} \mbox{and $\F(u_0,v_0) \le -K(m) \cdot (1+A^2)$} \bigg\},
  \end{eqnarray}
  the corresponding solution $(u,v)$ of \eqref{0} blows up at the finite time $\tm \in (0,\infty)$, where
  $\tm \le T(m,A)$.
\end{theo}

Actually, it turns out that when restricted to the power-type nonlinearities $\psi$ and $\phi$ and a uniformly parabolic system our result is optimal. Indeed, as is stated in the corollary below, in the case of uniformly parabolic system the infinite-time blowup does not hold. The finite-time blowup appears for the nonlinearities which are complementary to those leading to the boundedness of solutions proved in \cite{taowin_jde}.
\begin{cor}\label{cor1}
  Assume that $\phi(s) = (s+1)^{-p}$ and $\psi(s) =s (s+1)^{q-1}$, $s \ge 0$, with $p \le 0$ and $q \in \mathbb{R}$ such that
  $p+q > \frac{2}{n}$. Moreover, let
  $\Omega=B_R \subset \R^n$ with some $n\ge 3$ and $R>0$,  and let $m>0$ and $A>0$ be given.
  Then there exist positive constants $T(m,A)$ and $K(m)$ such that for any
  $(u_0,v_0) \in \B (m,A)$
  the corresponding solution $(u,v)$ of \eqref{0} blows up at the finite time $\tm \le T(m,A)$.
\end{cor}
Let us also mention that Corollary \ref{cor1} strenghtens the result in \cite{horstmann_winkler}. It shows that in a supercritical regime in \cite{horstmann_winkler} one really meets finite-time explosions of radially symmetric solutions. On the other hand, the situation considered in \cite{Miz_Sou}, where the blowup behavior of solutions is studied, is not covered by our result. This is a case of a critical exponent. 

Concerning the inifinite-time blowup, we have the following result.
\begin{theo}\label{theo2}
Let $\Omega \subset \R^n$ be a bounded domain with a smooth boundary, $n\geq 2$. Assume that there are $p,q \in \R$ with $p+2q < \frac{2}{n}$
and $D_1, D_2 >0$ such that for any $s\ge0$
\begin{equation}\label{balance}
  \phi (s) \ge D_1 (s+1)^{-p} \quad\mbox{and}\quad \psi (s) \le D_2 (s+1)^q
\end{equation}
are satisfied.
Then there exists a unique global-in-time solution $(u,v)$ to \eqref{0}.

Furthermore, if additionally $\Omega$ is a ball, \eqref{G1} and \eqref{H1} are fulfilled in the case $n \ge 3$ and
\begin{equation}\label{GH2}
  G(s) \le a s(\ln s)^\mu   \quad\mbox{and}\quad H(s) \le b \frac{s}{\ln s}, \; s \ge s_0,
\end{equation}
hold with some $s_0 >1$, $a,b>0$ and $\mu \in (0,1)$ in the case $n=2$, then there is a global-in-time radially symmetric solution $(u,v)$ to \eqref{0} which blows up in infinite time with respect to the $L^\infty(\Omega)$ norm.
\end{theo}

By the following corollary, which summarizes the results of Theorems \ref{theo1} and \ref{theo2}, we see that restricting ourselves to the power-type nonlinearities, there is still a gap to be filled when looking for critical exponents distinguishing between finite- and infinite-time blowup, namely $2/n>q\geq 0$.

\begin{cor}\label{cor2}
  Let $n \ge 2$. Assume that $\phi(s) = (s+1)^{-p}$ and $\psi(s) =s (s+1)^{q-1}$, $s \ge 0$, with $p,q \in \R$ and
  such that
  $p < \frac{2}{n} -2q$ is fulfilled. Then for $\Omega$ being a bounded domain with smooth boundary
  there exists a unique global-in-time solution.

  Moreover, let $\Omega=B_R \subset \R^n$, $R>0$. In case of $q<0$ and $\frac{2}{n} -q <p < \frac{2}{n}-2q$
  there exists initial data generating
  global-in-time radially symmetric solution $(u,v)$  to \eqref{0} which blows up at infinity with
  respect to the
  $L^\infty(\Omega)$ norm. On the other hand for $q>\frac{2}{n}$ and $p$ satisfying $p >\frac{2}{n}-q$ and, in the case 
	$q<1$ also $p \le 0$, there exists a
  radially symmetric solution $(u,v)$  to \eqref{0} which blows up at a finite time.
\end{cor}
Next, by \cite{win_mmas} we know that, taken initial data like in the introduction, there exists a classical unique
solution $(u,v)$ to \eqref{0}. A bound on the $L^\infty$ norm of $u$ on $(0,T)$ lets us prolong the classical solution to exist on $(0,\tm)$. Moreover, mass conservation for $u$ holds, as well as the bound $\int_\Omega vdx\leq A$.

In the following two sections we give the proofs of Theorems \ref{theo1} and \ref{theo2}. All the corollaries are simple applications of the theorems in the particular cases of power-type nonlinearities.

The last section is devoted to the application of our results in the studies of the so-called volume filling models, see \cite{hil:volume}, with a probability jump function of the form $q(u)=(1+u)^{-\gamma}$. We emphasize that in particular the existence of a critical mass phenomenon in dimension $2$ is shown. Namely, we indicate a value of mass $m_*$ such that solutions with initial mass smaller than $m_*$ are bounded. At the same time we point radially symmetric solutions with initial mass larger than $m_*$ that are unbounded. However, application of Corollary \ref{cor2} yields that such solutions exist for any time $t>0$. A critical mass phenomenon is important from the biological applications point of view, since it is often interpreted as a criterion of self-organization of cells.
\mysection{Finite-time blowup}\label{section3}

Our Theorem \ref{theo1} is an extension of our previous result in \cite{Cie_Sti_JDE}. It follows the strategy introduced in \cite{win_bu}. Here we give the sketch of the proof emphasizing only the steps which differ with respect
to \cite[Theorem 1.1]{Cie_Sti_JDE}.

The strategy of proving finite-time blowup is related to the Liapunov functional ${\cal F}$ associated to \eqref{0}. The required estimate which leads to a finite-time explosion of a solution is an inequality of the form
\begin{equation}\label{FF}
\frac{d}{dt}\left(-{\cal F} (u(t), v(t)) \right)\geq \Big(c(-{\cal F}) (u(t), v(t))-1\Big)^{\frac{1}{\theta}}
\end{equation}
for $t>0$ with some $\theta \in (0,1)$ and $c>0$. Indeed, once we know (\ref{FF}) we see that
\[
{\cal F}(u(t),v(t))\rightarrow -\infty
\]
as $t\rightarrow \bar{T}$, for some $\bar{T}<\infty$  provided the initial value $-{\cal F} (u_0, v_0)$ is large enough. But once we know that ${\cal F}$ tends to $-\infty$ at a finite time, we are sure that $\int_\Omega uv dx$ tends to $\infty$ as this integral is the only negative ingredient of ${\cal F}$ (see (\ref{F}) and (\ref{GH})). Since unboundedness of $\int_\Omega uv dx$ along with the boundedness of $\Omega$ yields finite-time blowup of either $u$ or $v$ in $L^\infty$, $u$ blows up in finite time.

In order to describe the strategy precisely we introduce the following notation. We fix $m>0$, $M>0$, $B>0$, and $\kappa>n-2$ and
assume that
\be{m}
	\io u = m \qquad \mbox{and} \qquad \io v \le M
\ee
and
\be{B}
	v(x) \le B|x|^{-\kappa} \qquad \mbox{for all } x\in\Omega
\ee
are fulfilled. Furthermore, we define the space
\bea{S}
	\set &:=& \bigg\{ (u,v) \in C^1(\bar\Omega) \times C^2(\bar\Omega) \ \bigg| \
	\mbox{$u$ and $v$ are positive and radially} \nn\\
	& & \hspace*{5mm}
	\mbox{symmetric satisfying $\frac{\partial v}{\partial\nu}=0$ on $\pO$, \eqref{m}, and \eqref{B}} \bigg\}.
\eea
Next, for given $f$ and $g$, we consider the system
\be{f}
-\Delta v + v - u=f,
\ee
\be{g}
\left( \frac{\phi (u)}{\sqrt{\psi(u)}}\nabla u -\sqrt{\psi (u)}\nabla v \right) \cdot \frac{x}{|x|}=g, \qquad x \neq 0.
\ee
We notice that proving the inequality
  \begin{eqnarray}\label{crucial}
	\io uv &\le&  C(m,\kappa) \cdot \left( 1+M^2 + B^{\frac{2n+4}{n+4}} \right) \cdot
	\Bigg( \|f\|_{L^2(\Omega)}^{2\theta}
	+ \left\|g\right\|_{L^2(\Omega)}^{2\theta} +1 \Bigg)\nn \\
    &+& \frac{1}{2} \io | \nabla v |^2 + \io G(u),
  \end{eqnarray}
with
\be{theta}
	\theta:=\frac{1}{1+\frac{n}{(2n+4)\kappa}} \,
	\in \Big(\frac{1}{2},1\Big)
  \ee
for all solutions $(u,v)\in \set$ of the hyperbolic-elliptic system (\ref{f})-(\ref{g}), implies that inequality
\begin{equation}\label{4.1}
  \frac{\F(u,v)}{\D^\theta(u,v)+1} \ge - C(m,M,B) \qquad \mbox{for all } (u,v) \in \set
\end{equation}
is satisfied with some constants $\theta \in (0,1)$ and $C(m,M,B)>0$ . In view of (\ref{liapunov}) and the fact that we can choose properly initial data such that $-{\cal F}(u_0,v_0)$ is large enough, we arrive at (\ref{FF}) and the finite-time blowup result is proved provided solutions to \eqref{0} belong to $\set$. Both, the possibility of choosing $-{\cal F}(u_0,v_0)$ large enough and the fact that solutions to \eqref{0} belong to $\set$, were already proved in \cite{Cie_Sti_JDE}. Hence we are left with the proof of (\ref{crucial}) in the set $\set$. This is done in a few steps, where we follow the corresponding proof in \cite{Cie_Sti_JDE}.

We observe that for $\eps>0$ (see \cite[Lemma 3.2]{Cie_Sti_JDE})
\be{3.2.1}
	\io uv \le (1+ \eps) \io |\nabla v|^2 + C(\eps) \cdot \left( 1+ M^2 \right) \cdot \left(
	\Big\| f\Big\|_{L^2(\Omega)}^{\frac{2n+4}{n+4}} +1 \right).
  \ee

Next, we modify further steps in the proof of \cite[Lemma 3.1]{Cie_Sti_JDE}. First we prove the lemma
corresponding to \cite[Lemma~3.4]{Cie_Sti_JDE}. Here we essentially use \eqref{psi} instead of the superlinear growth of
$\psi$ which was assumed in \cite{Cie_Sti_JDE}.
\begin{lem}\label{lem3.4}
  Assume that \eqref{H1} and \eqref{psi} are satisfied. Then there exist $\mu = \mu(\gamma) \in (0,2)$
  and $C(m)>0$ such that for all $r_0\in (0,R)$ and $(u,v)\in\set$
  \begin{eqnarray}\label{3.4.1}
	\int_{B_{r_0}} |\nabla v|^2
	&\le& \mu \io G(u) + C(m) \cdot \Bigg\{
	r_0 \cdot \Big\|\Delta v-v+u \Big\|_{L^2(\Omega)}^2 \nn \\
	& & + r_0 \cdot \left\| \frac{\phi (u)}{\sqrt{\psi(u)}}\nabla u -\sqrt{\psi (u)}\nabla v \right\|_{L^2(\Omega)}^2
	+\|v\|_{L^2(\Omega)}^2
	+1 \Bigg\}
  \end{eqnarray}
  is fulfilled.
\end{lem}
\proof As \eqref{H1} implies $(\frac{4(n-1)}{n-2} -2)\gamma <2$, we can fix $\delta \in (0, \frac{2n-2}{R}]$ small
enough such that
\begin{equation}\label{mu}
 \mu_1 := \left( \frac{4(n-1)}{n-2} e^{\delta R} -2 \right) \cdot \gamma  \in (0,2)
\end{equation}
is fulfilled. Next we fix $\eta := \frac{2-\mu_1}{2R}$ and observe that
\begin{equation}\label{2.12}
\mu := \mu_1 + \eta R \in (0,2)
\end{equation} is fulfilled.

Our proof follows the lines of the proof of \cite[Lemma 3.4]{Cie_Sti_JDE}. After rewriting (\ref{f})-(\ref{g}) in radial coordinates we arrive at the estimate \cite[(3.25)]{Cie_Sti_JDE}:
\begin{eqnarray}\label{3.4.3}
 r^{2n-2} v_r^2(r) &\le&
	-2\int_0^r e^{\delta (r-\rho)} \rho^{2n-2} \frac{u(\rho) \phi(u(\rho))}{\psi(u(\rho))} u_r (\rho) 
        \; d\rho \nn \\ & &
	+2\int_0^r e^{\delta (r-\rho)} \rho^{2n-2} \frac{u(\rho)}{\sqrt{\psi(u(\rho))}} g (\rho) \; d\rho \nn\\
	& & +\frac{1}{\delta} \int_0^r e^{\delta (r-\rho)} \rho^{2n-2} f^2(\rho) d\rho
	+\int_0^r e^{\delta (r-\rho)} \rho^{2n-2} (v^2)_r(\rho) d\rho, \; r \in (0,R).
\end{eqnarray}
Out of four terms on the right-hand side of \eqref{3.4.3} only
%\[
%2\int_0^r e^{\delta (r-\rho)} \rho^{2n-2} \frac{u(\rho)}{\sqrt{\psi(u(\rho))}} g (\rho) \; d\rho
%\]
the second one requires a more detailed treatment than in \cite{Cie_Sti_JDE}.  Denoting by $\omega_n$ the $(n-1)$-dimensional measure of the sphere $\partial B_1$ and applying the
Cauchy-Schwarz inequality as well as \eqref{psi} and Young's inequality, we arrive at
\begin{eqnarray}\label{3.4.5}
 && \hspace*{-20mm}
 2\int_0^r e^{\delta (r-\rho)} \rho^{2n-2} \frac{u(\rho)}{\sqrt{\psi(u(\rho))}} g (\rho) \; d\rho \nn \\
 &\le& 2 \left( \int_0^R \rho^{n-1} \frac{u^2(\rho)}{\psi(u(\rho))} \; d\rho \right)^\frac{1}{2} \cdot
	\left( \int_0^r e^{2\delta(r-\rho)} \cdot \rho^{3n-3} g^2(\rho) \; d\rho \right)^\frac{1}{2} \nn \\
 &\le& 2 \left( L \int_0^R \rho^{n-1} \left( G(u(\rho)) + u(\rho) +1 \right) \; d\rho \right)^\frac{1}{2} \cdot
	\left( e^{2\delta R} r^{2n-2} \int_0^R \rho^{n-1} g^2(\rho) \; d\rho \right)^\frac{1}{2} \nn \\
 &\le& \frac{2 \sqrt{L} \; e^{\delta R}}{w_n} \cdot \left( \io (G(u) +u+1) \right)^{\frac{1}{2}} \cdot r^{n-1} \cdot \| g \|_{L^2(\Omega)}
 \nn \\
 &\le& \frac{\eta}{\omega_n} r^{n-1} \left( m+|\Omega| + \io G(u) \right) + \frac{L e^{2\delta R}}{\eta w_n}  r^{n-1} \cdot \| g \|_{L^2(\Omega)}^2,
 \qquad r \in (0,R).
\end{eqnarray}
Hence, estimating the other terms on the right-hand side of \cite[(3.25)]{Cie_Sti_JDE} like in \cite[(3.26), (3.28) and (3.29)]{Cie_Sti_JDE} and using (\ref{3.4.5}), we see that there is a constant $c_1(m) >0$ such that
\begin{eqnarray*}
  r^{2n-2} v_r^2(r) &\le& 4(n-1) e^{\delta R} \int_0^r \rho^{2n-3} H(u(\rho))  \; d\rho -2r^{2n-2} H(u(r))
	+ \frac{\eta}{\omega_n} r^{n-1} \io G(u) \\
  & & + \frac{c_1(m)}{\omega_n} r^{n-1} + \frac{c_1(m)}{\omega_n} r^{n-1} \|g\|_{L^2(\Omega)}^2
   + \frac{c_1(m)}{\omega_n} r^{n-1} \|f\|_{L^2(\Omega)}^2
	+ r^{2n-2} v^2(r)
\end{eqnarray*}
for $r \in (0,R)$. Multiplying this inequality by $\omega_n r^{1-n}$ and integrating with respect to $r \in (0,r_0)$, we have
\begin{eqnarray}\label{3.4.8}
 \int_{B_{r_0}} | \nabla v |^2 &=& \omega_n \int_0^{r_0} r^{n-1} v_r^2(r) \; dr \nn \\
 &\le& 4(n-1) e^{\delta R} \omega_n  \int_0^{r_0} r^{1-n} \int_0^r \rho^{2n-3} H(u(\rho))  \; d\rho \, dr \nn \\
 & &  -2\omega_n \int_0^{r_0} r^{n-1} H(u(r)) \; dr + \eta r_0 \io G(u) + c_1(m) r_0 \nn \\ 
 & & + c_1 (m) r_0 \|g\|_{L^2(\Omega)}^2  + c_1(m) r_0 \|f\|_{L^2(\Omega)}^2 + \omega_n \int_0^{r_0} r^{n-1} v^2(r) \; dr \nn \\
 &\le& 4(n-1) e^{\delta R} \omega_n  \int_0^{r_0} r^{1-n} \int_0^r \rho^{2n-3} H(u(\rho))  \; d\rho \, dr
 -2 \int_{B_{r_0}} H(u) \nn \\
 & &  + \eta R \io G(u) + c_1 (m) \left( R +  r_0 \|g\|_{L^2(\Omega)}^2
 + r_0 \|f\|_{L^2(\Omega)}^2 \right) + \|v\|_{L^2(\Omega)}^2.
\end{eqnarray}\
Finally, Fubini's theorem, $n \ge 3$, the nonnegativity of $H$ and \eqref{H1} yield
\begin{eqnarray}\label{3.4.81}
 && \hspace*{-20mm}
   4(n-1) e^{\delta R} \omega_n  \int_0^{r_0} r^{1-n} \int_0^r \rho^{2n-3} H(u(\rho))  \; d\rho \, dr
   -2 \int_{B_{r_0}} H(u) \nn\\
 &=& 4(n-1) e^{\delta R} \omega_n  \int_0^{r_0} \left( \int_{\rho}^{r_0} r^{1-n} \; dr \right)
    \rho^{2n-3} H(u(\rho))  \; d\rho -2 \int_{B_{r_0}} H(u) \nn \\
 &=& \frac{4(n-1)}{n-2} e^{\delta R} \omega_n  \int_0^{r_0} \left( \rho^{2-n} - r_0^{2-n} \right)
    \rho^{2n-3} H(u(\rho))  \; d\rho -2 \int_{B_{r_0}} H(u) \nn \\
 &\le& \frac{4(n-1)}{n-2} e^{\delta R} \omega_n  \int_0^{r_0}
    \rho^{n-1} H(u(\rho))  \; d\rho -2 \int_{B_{r_0}} H(u) \nn \\
 &=& \left( \frac{4(n-1)}{n-2} e^{\delta R} -2 \right) \int_{B_{r_0}} H(u)
 \le \left( \frac{4(n-1)}{n-2} e^{\delta R} -2 \right) \int_{\Omega} H(u) \nn \\
 &\le& \left( \frac{4(n-1)}{n-2} e^{\delta R} -2 \right) \int_{\Omega} \left( \gamma G(u) + b(u+1) \right)
 = \mu_1 \io G(u) + c_2 (m),
\end{eqnarray}
where $\mu_1$ was defined in (\ref{mu}). Hence, a combination of (\ref{3.4.8}) and (\ref{3.4.81}) yields
\begin{eqnarray}
\int_{B_{r_0}} | \nabla v |^2 &\le& \mu_1 \io G(u) + c_2 (m) + \eta R \io G(u) + c_1(m) R \nn \\
 & & + c_1 (m) r_0 \|g\|_{L^2(\Omega)}^2
 + c_1(m) r_0 \|f\|_{L^2(\Omega)}^2 + \|v\|_{L^2(\Omega)}^2.
\end{eqnarray}
In view of (\ref{2.12}), the claim \eqref{3.4.1} is proved.
\qed

The next lemma corresponds to \cite[Lemma~3.5]{Cie_Sti_JDE}.
\begin{lem}\label{lem3.5}
  Suppose that \eqref{H1} and \eqref{psi} are fulfilled and let $\theta \in (\frac{1}{2},1)$ and
  $\mu \in (0,2)$ be as defined in \eqref{theta} and Lemma~\ref{lem3.4},
  respectively. Then for any $\eps \in (0, \frac{1}{2})$ there exists $C(\eps,m,\kappa)>0$ such that
  \begin{eqnarray}\label{3.5.1}
	\io |\nabla v|^2
	&\le&  C(\eps,m,\kappa) \cdot \left( 1+ M^2
    + B^{\frac{2n+4}{n+4}} \right) \cdot
	\bigg( \Big\|\Delta v-v+u\Big\|_{L^2(\Omega)}^{2\theta} \nn \\ & &
	+ \Big\|\frac{\phi (u)}{\sqrt{\psi(u)}}\nabla u -\sqrt{\psi (u)}\nabla v\Big\|_{L^2(\Omega)}^{2\theta} +1 \bigg) \nn \\ &&
	+ \frac{\eps}{1-2\eps} \io uv + \frac{\mu}{1-2\eps} \io G(u)
  \end{eqnarray}
  is fulfilled for all $(u,v) \in \set$.
\end{lem}
\proof We fix $\eps \in (0, \frac{1}{2})$ and set $\beta := \frac{(2n+4)\kappa}{n}$ which implies
$\theta = \frac{\beta}{\beta+1}$. Next we define $r_0 :=
\min\{\frac{R}{2}, \left( \|f\|_{L^2(\Omega)} + \|g\|_{L^2(\Omega)}\right)^{-\frac{2}{\beta+1}} \} \in (0,R)$. Hence, by
\cite[Lemma~3.3]{Cie_Sti_JDE} there
is $c_1= C_1 (\eps, m, \kappa) \cdot \big( 1+ M^{\frac{2n+4}{n+4}} + B^{\frac{2n+4}{n+4}} \big)>0$ such that
\begin{equation}\label{3.5.2}
 \int_{\Omega \setminus B_{r_0}} |\nabla v|^2
	\le \eps \io uv + \eps \io |\nabla v|^2
	+ c_1 \cdot \Big( r_0^{-\beta} + \|f\|_{L^2(\Omega)}^\frac{2n+4}{n+4} \Big).
\end{equation}
Applying next Lemma~\ref{lem3.4}, we get a constant $c_2 = c_2(m)$ such that
\begin{equation}\label{3.5.3}
 \int_{B_{r_0}} |\nabla v|^2  \le \mu \io G(u) + c_2 \cdot \Big( r_0 \|f\|_{L^2(\Omega)}^2 + r_0 \|g\|_{L^2(\Omega)}^2
   + \|v\|_{L^2(\Omega)}^2 + 1 \Big).
\end{equation}
Adding both inequalities, we deduce that
\begin{eqnarray}\label{3.5.4}
 (1-\eps) \io |\nabla v|^2  &\le& \eps \io uv + \mu \io G(u) + c_1 r_0^{-\beta} +
  c_1 \|f\|_{L^2(\Omega)}^\frac{2n+4}{n+4}  \nn \\
	& & + c_2 r_0 \left( \|f\|_{L^2(\Omega)}^2 + \|g\|_{L^2(\Omega)}^2 \right)
  + c_2 + c_2 \|v\|_{L^2(\Omega)}^2.
\end{eqnarray}
Next, by \cite[Lemma~2.2]{Cie_Sti_JDE} and \eqref{m} there exists $c_3 = C_3 (\eps,m) \cdot M^2 >0$ such that
$$c_2 \|v\|_{L^2(\Omega)}^2 \le \eps \io |\nabla v|^2 + c_3,$$
which inserted into \eqref{3.5.4} yields
\begin{eqnarray}\label{3.5.5}
 (1-2\eps) \io |\nabla v|^2  &\le& \eps \io uv + \mu \io G(u) + c_2 + c_3  +I,
\end{eqnarray}
where we set
$$I := c_1 r_0^{-\beta} +
  c_1 \|f\|_{L^2(\Omega)}^\frac{2n+4}{n+4}  + c_2 r_0 \left( \|f\|_{L^2(\Omega)}^2 + \|g\|_{L^2(\Omega)}^2 \right).$$
In case of $\|f\|_{L^2(\Omega)} + \|g\|_{L^2(\Omega)} \le (\frac{2}{R})^\frac{\beta+1}{2}$, we have $r_0=\frac{R}{2}$ and
conclude that
  \bas
	I \le c_1 \cdot \Big(\frac{2}{R}\Big)^\beta
	+ c_1 \cdot \Big(\frac{2}{R}\Big)^{\frac{\beta+1}{2} \cdot \frac{2n+4}{n+4}}
	+ c_2 \cdot \frac{R}{2} \cdot 2 \Big( \frac{2}{R} \Big)^{\beta+1},
  \eas
  which in conjunction with \eqref{3.5.5} proves \eqref{3.5.1} in this case.

Furthermore, in the case $\|f\|_{L^2(\Omega)} + \|g\|_{L^2(\Omega)} > (\frac{2}{R})^\frac{\beta+1}{2}$ we have
$r_0= (\|f\|_{L^2(\Omega)} + \|g\|_{L^2(\Omega)})^{-\frac{2}{\beta+1}}$ and therefore
\begin{eqnarray*}
  I &\le& c_1 \left( \|f\|_{L^2(\Omega)} + \|g\|_{L^2(\Omega)} \right)^\frac{2\beta}{\beta+1}
	+c_1 \|f\|_{L^2(\Omega)}^\frac{2n+4}{n+4}
	+ c_2 \|f\|_{L^2(\Omega)}^{2-\frac{2}{\beta+1}} + c_2 \|g\|_{L^2(\Omega)}^{2-\frac{2}{\beta+1}} \\
	&=& (4c_1+c_2) \|f\|_{L^2(\Omega)}^\frac{2\beta}{\beta+1}
	+c_1 \|f\|_{L^2(\Omega)}^\frac{2n+4}{n+4} + (4c_1+c_2) \|g\|_{L^2(\Omega)}^\frac{2\beta}{\beta+1}.
\end{eqnarray*}	
In view of $\kappa > n-2$ and $n \ge 3$, we calculate
$$\frac{\beta}{\frac{n+2}{2}} =\frac{2}{n+2} \cdot \frac{(2n+4)\kappa}{n}
	>\frac{4(n-2)}{n} \ge \frac{4}{3}>1$$
which implies that $2\theta = \frac{2\beta}{\beta+1}>\frac{2n+4}{n+4}$. Applying once more Young's inequality,
we obtain
$$I \le (5c_1+c_2) \|f\|_{L^2(\Omega)}^\frac{2\beta}{\beta+1} +  c_1 + (4c_1+c_2) \|g\|_{L^2(\Omega)}^\frac{2\beta}{\beta+1},$$
which inserted into \eqref{3.5.5} proves \eqref{3.5.1} in the case
  $\|f\|_{L^2(\Omega)} + \|g \|_{L^2 (\Omega)}> (\frac{2}{R})^\frac{\beta+1}{2}$ and thereby completes the proof.
\qed

Finally, we complete the proof of the announced estimate \eqref{crucial} like in \cite[Proof of Lemma~3.1]{Cie_Sti_JDE} (just replacing
$\|g\|_{L^2(\Omega)}$ by $\|g\|_{L^2(\Omega)}^{2\theta}$). The proof of Theorem~\ref{theo1} is then completely the same as given in
Theorem~3.6 and Section~4 of \cite{Cie_Sti_JDE} for proving \cite[Theorem~1.1]{Cie_Sti_JDE}. \\

\mysection{Infinite-time blowup}\label{section4}
This section is devoted to the proof of Theorem~\ref{theo2}. %Energy estimates that we give are restricted to convex domains. 
The essential part of the energy estimates relies on a lemma, which holds in bounded domains $\Omega$ with smooth
boundary, stating that a function which satisfies the homogeneous Neumann boundary condition on $\partial\Omega$ fulfills
\begin{equation}\label{convex}
\frac{\partial |\nabla v|^2}{\partial \nu} \le c_\Omega |\nabla v|^2\;\;\mbox{on}\;\; \partial \Omega
\end{equation}
with some constant $c_\Omega >0$ depending only on the curvatures of $\partial \Omega$,
see for instance \cite[Lemma~4.2]{Miz_Sou}. This generalizes the estimate $\frac{\partial |\nabla v|^2}{\partial \nu} \le 0$ on
$\partial \Omega$ which holds for convex domains only and is known in thermodynamics, see for example \cite[Lemma on p.95]{evans}
or \cite[Lemma~3.2]{taowin_jde}. The first use of energy estimates relying on the latter estimate in chemotaxis goes back to \cite{winkler_cpde}. In what follows we use many ideas from the proof of \cite[Lemma~3.3]{taowin_jde}, which is the application of the above mentioned estimates to the quasilinear Keller-Segel system, and combine them with the ideas from \cite{ishida}
in order to remove the convexity assumption on $\Omega$ with the help of \eqref{convex}.
\begin{lem}\label{lemat}
Let $\Omega\subset \R^n$ be a bounded domain with smooth boundary for some $n\ge 2$. Moreover, assume that \eqref{balance} holds with
$p+2q <\frac{2}{n}$.
Then for any solution $(u,v)$ to \eqref{0}, any $\gamma \in [1,\infty)$ and any $T \in (0,
\infty)$ with $T \le \tm$ there is $C>0$ such that $u$ admits the estimate
\begin{equation}\label{wzor}
\|u(\cdot,t)\|_{L^\gamma(\Omega)}\leq C, \qquad t \in \left( 0,T \right).
\end{equation}
\end{lem}

\proof We fix $\gamma > 1$ and $\alpha >1$ which will be specified later. Then we multiply the first equation in \eqref{0} by
  $(u+1)^{\gamma-1}$ and use \eqref{balance} as well as Young's inequality to obtain
  \begin{eqnarray}\label{4.1.1}
	  & & \hspace*{-15mm} \frac{d}{dt} \io (u+1)^\gamma \nn \\
		&=& \gamma \io (u+1)^{\gamma -1} \nabla \cdot ( \phi(u) \nabla u - \psi(u) \nabla v) \nn \\
		&=& - \gamma (\gamma-1) \io (u+1)^{\gamma -2} \phi(u) |\nabla u|^2 + \gamma (\gamma -1) \io (u+1)^{\gamma-2} \psi(u) \nabla u
		\nabla v \nn \\
		&\le& - D_1 \gamma (\gamma-1) \io (u+1)^{\gamma -2-p}|\nabla u|^2 + D_2\gamma (\gamma -1) \io (u+1)^{\gamma-2+q} |\nabla u
		\nabla v| \nn \\
		&\le& - \frac{D_1 \gamma (\gamma-1)}{2} \io (u+1)^{\gamma -2-p}|\nabla u|^2
		+ \frac{D_2^2\gamma (\gamma -1)}{2D_1} \io (u+1)^{\gamma + p+2q-2} |\nabla v|^2.
	\end{eqnarray}
In view of the second equation in \eqref{0} and the identity
\[
\Delta|\nabla v|^2-2|D^2v|^2=2\nabla v\cdot \nabla \Delta v
\]
we have
	\begin{eqnarray}\label{4.1.2}
	  \frac{d}{dt} \io |\nabla v|^{2\alpha}
		&=& \alpha \io |\nabla v|^{2\alpha -2} 2 \nabla v \nabla v_t
		=  \alpha \io |\nabla v|^{2\alpha -2} \left( 2 \nabla v \nabla \Delta v - 2 |\nabla v|^2 + 2 \nabla u \nabla v \right) \nn \\
		&=& \alpha \io |\nabla v|^{2\alpha -2} \left( \Delta |\nabla v|^2 - 2 |D^2 v|^2 - 2 |\nabla v|^2 + 2 \nabla u \nabla v \right)
		\nn \\
		&\le& - \alpha (\alpha-1) \io |\nabla v|^{2\alpha -4} \left| \nabla |\nabla v|^2 \right|^2
		+ \alpha c_\Omega \int_{\partial \Omega} |\nabla v|^{2\alpha} d \sigma \nn \\ 
		& & - 2\alpha \io |\nabla v|^{2\alpha -2} |D^2 v|^2 - 2\alpha \io |\nabla v|^{2\alpha}
		+ 2\alpha \io |\nabla v|^{2\alpha -2} \nabla u \nabla v,
	\end{eqnarray}
where the last inequality is obtained by integration by parts and the use of (\ref{convex}). Next an integration by parts in conjunction with Young's inequality and the estimate $|\Delta v|^2 \le n |D^2 v|^2$
	show that
	\begin{eqnarray*}
	  2\alpha \io |\nabla v|^{2\alpha -2} \nabla u \nabla v
		&=& - 2\alpha (\alpha-1) \io u |\nabla v|^{2\alpha -4} \nabla v \nabla |\nabla v|^2 - 2\alpha \io u |\nabla v|^{2\alpha-2}
		\Delta v \\
		&\le& \frac{\alpha (\alpha-1)}{2} \io |\nabla v|^{2\alpha -4} \left| \nabla |\nabla v|^2 \right|^2
		+ 2 \alpha (\alpha-1) \io u^2 |\nabla v|^{2\alpha -2} \\ & &  + \frac{2\alpha}{n} \io |\nabla v|^{2\alpha-2} |\Delta v|^2
		+ \frac{\alpha n}{2} \io u^2 |\nabla v|^{2\alpha -2} \\
		&\le& \frac{\alpha (\alpha-1)}{2} \io |\nabla v|^{2\alpha -4} \left| \nabla |\nabla v|^2 \right|^2
		+ 2\alpha \io |\nabla v|^{2\alpha-2} |D^2 v|^2 \\
		& & + \left( 2 \alpha (\alpha-1) + \frac{\alpha n}{2} \right) \io u^2 |\nabla v|^{2\alpha -2}.
	\end{eqnarray*}
	Inserting this into \eqref{4.1.2}, we deduce that
	\begin{eqnarray}\label{4.1.3}
	  \frac{d}{dt} \io |\nabla v|^{2\alpha} &\le& -\frac{2(\alpha-1)}{\alpha} \io \left| \nabla |\nabla v|^\alpha \right|^2
		+ \left( 2 \alpha (\alpha-1) + \frac{\alpha n}{2} \right) \io u^2 |\nabla v|^{2\alpha -2} \nn \\
		& & + \alpha c_\Omega \int_{\partial \Omega} |\nabla v|^{2\alpha} d \sigma.
	\end{eqnarray}
	Adding \eqref{4.1.1} and \eqref{4.1.3} and using $p+2q <2$ in conjunction with Young's inequality, we thus conclude that
	\begin{eqnarray}\label{4.1.4}
	  & & \hspace*{-15mm} \frac{d}{dt} \left( \io (u+1)^\gamma + \io |\nabla v|^{2\alpha} \right)
		+ \frac{2(\alpha-1)}{\alpha} \io \left| \nabla |\nabla v|^\alpha \right|^2 \nn \\
		&\le& c_1 \io (u+1)^{\gamma + p +2q-2} |\nabla v|^2 + c_2 \io (u+1)^2 |\nabla v|^{2\alpha -2} 
		+ \alpha c_\Omega \int_{\partial \Omega} |\nabla v|^{2\alpha} d \sigma \nn \\
		&\le& c_3 \io (u+1)^\gamma + c_4 \io |\nabla v|^{\frac{2\gamma}{2-p-2q}} + c_5 \io |\nabla v|^{\frac{2(\alpha -1) \gamma}
		{\gamma -2}} + \alpha c_\Omega \int_{\partial \Omega} |\nabla v|^{2\alpha} d \sigma,
	\end{eqnarray}
	if we choose $\gamma > \gamma_0 := \max \{ 2, 2-p-2q\}$. Furthermore, as $\| u(t) \|_{L^1 (\Omega)} = \| u_0 \|_{L^1 (\Omega)}$
	for any $t \in (0, \tm)$ due to the conservation of mass and $v_0 \in W^{1,\infty} (\Omega)$, we infer from
	\cite[Lemma~4.1]{horstmann_winkler} that for
	any $\delta \in (1, \frac{n}{n-1})$ there exists $C_\delta >0$ such that $\| \nabla v(t) \|_{L^\delta (\Omega)} \le C_\delta$
	for all $t \in (0, \tm)$. 
	Hence, we fix $r \in (0,\frac{1}{2})$ and see that for any $\alpha >1$ and 
	$s \in (0, \frac{n}{(n-1)\alpha})$ we have
	$$a_1 := \frac{\frac{1}{s} - \frac{1}{2}+ \frac{1}{2n} + \frac{r}{n}}{\frac{1}{n} -\frac{1}{2} + \frac{1}{s}} \in (0,1).$$
	Using as in \cite{ishida} the compact embedding of $W^{r+ \frac{1}{2},2} (\Omega)$ into $L^2 (\partial \Omega)$ 
	(see \cite[Proposition~4.22(ii) and Theorem~4.24(i)]{ht}) and the fractional 
	Gagliardo-Nirenberg inequality (see \cite[Lemma~2.5]{ishida}), we obtain
	\begin{eqnarray}\label{4.1.4a}
	  \alpha c_\Omega \int_{\partial \Omega} |\nabla v|^{2\alpha} d \sigma &=&
		\alpha c_\Omega \| |\nabla v|^{\alpha} \|_{L^2 (\partial \Omega)}^2 
		\le C \| |\nabla v|^{\alpha} \|_{W^{r+ \frac{1}{2},2} (\Omega)}^2 \nn \\
		&\le& C \left( \| \nabla |\nabla v|^\alpha \|_{L^2 (\Omega)}^{2 a_1}
		\| |\nabla v|^\alpha \|_{L^s (\Omega)}^{2(1-a_1)}
		+ \| |\nabla v|^\alpha \|_{L^s (\Omega)}^2 \right) \nn \\
		&\le& C \left(1+ \| \nabla |\nabla v|^\alpha \|_{L^2 (\Omega)}^{2a_1} \right)
		\le \frac{\alpha-1}{2\alpha} \| \nabla |\nabla v|^\alpha \|_{L^2 (\Omega)}^2 + C.
	\end{eqnarray}
	Furthermore,	if we choose $\beta >0$ and $s \in (0, \frac{n}{(n-1)\alpha})$ such that
	\begin{equation}\label{4.1.5}
	  a := \frac{\frac{1}{s} - \frac{\alpha}{\beta}}{\frac{1}{n} - \frac{1}{2} + \frac{1}{s}} \in (0,1)
		\quad\mbox{and}\quad \frac{a\beta}{\alpha} <2
	\end{equation}
	the Gagliardo-Nirenberg inequality implies
	\begin{eqnarray}\label{4.1.6}
	  \io |\nabla v|^\beta &=& \| |\nabla v|^\alpha \|_{L^{\frac{\beta}{\alpha}} (\Omega)}^{\frac{\beta}{\alpha}} \nn \\
		&\le& C \left( \| \nabla |\nabla v|^\alpha \|_{L^2 (\Omega)}^{\frac{a\beta}{\alpha}}
		\| |\nabla v|^\alpha \|_{L^s (\Omega)}^{\frac{(1-a)\beta}{\alpha}}
		+ \| |\nabla v|^\alpha \|_{L^s (\Omega)}^{\frac{\beta}{\alpha}} \right) \nn \\
		&\le& C \left(1+ \| \nabla |\nabla v|^\alpha \|_{L^2 (\Omega)}^{\frac{a\beta}{\alpha}} \right)
		\le \frac{\alpha-1}{2\alpha} \| \nabla |\nabla v|^\alpha \|_{L^2 (\Omega)}^2 + C.
	\end{eqnarray}
	In view of \eqref{4.1.4} we would like to use \eqref{4.1.6} for $\beta_1 := \frac{2\gamma}{2-p-2q}$ and $\beta_2 :=
	\frac{2(\alpha -1) \gamma}{\gamma -2}$. Hence we have to ensure that we can choose $\alpha$, $\gamma$ and $s$ appropriately
	such that \eqref{4.1.5} is satisfied for $\beta_1$ and $\beta_2$.
	
To this end let us define
\[
a(\alpha, \gamma, s):=\frac{\frac{1}{s}-\frac{\alpha}{\beta(\alpha, \gamma)}}{\frac{1}{n}-\frac{1}{2}+\frac{1}{s}}
\]
and
\[
f(\alpha, \gamma, s):=\frac{a(\alpha, \gamma, s)\beta(\alpha, \gamma)}{\alpha}.
\]
We first observe that when $\alpha$ and $\gamma$ are fixed, $a$ and $f$ are continuous in $s$ in a neighborhood of   $s = \frac{n}{(n-1)\alpha}$. Hence, once we can identify $\alpha'$ and $\gamma'$ such that (\ref{4.1.5}) holds with
$a(\alpha', \gamma',\frac{n}{(n-1)\alpha})$ and $f(\alpha', \gamma',\frac{n}{(n-1)\alpha})$, by continuity, keeping the same values $\alpha'$ and $\gamma'$, we can pick $s \in (0, \frac{n}{(n-1)\alpha})$ close to $\frac{n}{(n-1)\alpha}$ such that \eqref{4.1.5} holds for $\beta_1$ and $\beta_2$.

Thus, it is enough to focus on the case $s = \frac{n}{(n-1)\alpha}$. Then, $a \in (0,1)$ and $\frac{a\beta}{\alpha} <2$ if
	\begin{equation}\label{4.1.7}
	  \beta > \frac{n}{n-1} \quad\mbox{and} \quad \alpha > \frac{\beta}{2} - \frac{1}{n-1}.
	\end{equation}
Indeed, in view of $\alpha >1$ and $n \ge 2$, the first inequality in (\ref{4.1.7}) is equivalent to $a(\alpha, \gamma,\frac{n}{(n-1)\alpha})>0$.
Next, we notice that $a(\alpha, \gamma,\frac{n}{(n-1)\alpha})<1$ is equivalent to
\[
\frac{\alpha}{\beta}>\frac{1}{2}-\frac{1}{n},
\]
and the latter estimate is a consequence of (\ref{4.1.7}). Moreover, $f(\alpha, \gamma,\frac{n}{(n-1)\alpha})<2$ holds if
\[
\frac{\beta}{2} \left(1-\frac{1}{n}-\frac{1}{\beta}\right)<\frac{1}{n}-\frac{1}{2}+\alpha\left(1-\frac{1}{n}\right)
\]
is satisfied,	while the latter is equivalent to the second inequality in \eqref{4.1.7}.
	
Now \eqref{4.1.7} is satisfied for $\beta_1$ and $\beta_2$ if
	\begin{equation}\label{4.1.8}
	  \gamma > \frac{(2-p-2q)n}{2(n-1)}, \quad \alpha > \frac{n}{2(n-1)} +1,
	\end{equation}
and	
	\begin{equation}\label{4.1.81}	
		\gamma \frac{1}{2-p-2q} -\frac{1}{n-1} < \alpha < \gamma \frac{n}{2(n-1)} - \frac{1}{n-1}
	\end{equation}
are fulfilled. Indeed, the inequalities in (\ref{4.1.8}) guarantee the first inequality in \eqref{4.1.7} for $\beta_1$ and $\beta_2$, respectively. On the other hand, the left inequality in \eqref{4.1.81} is simply the second inequality in \eqref{4.1.7} for $\beta_1$, while the right inequality in \eqref{4.1.81} is equivalent to the second inequality in \eqref{4.1.7} for $\beta_2$.

Hence, we first fix $\gamma_1 := \max \{ \gamma_0, \frac{3(2-p-2q)n}{2(n-1)} \}$. Then, for an arbitrary $\gamma \ge \gamma_1$ we choose
$\alpha$ such that \eqref{4.1.81} is satisfied, where the latter is possible as $p+2q <\frac{2}{n}$ implies $\frac{1}{2-p-2q} <
\frac{n}{2(n-1)}$. Moreover, $\gamma$ and $\alpha$ satisfy \eqref{4.1.8} in view of our choice of $\gamma_1$.
Hence, we conclude that for any $\gamma \ge \gamma_1$ there exists $\alpha >1$ such that \eqref{4.1.8} and \eqref{4.1.81} are fulfilled.

Therefore, using \eqref{4.1.6} for $\beta_1$ and $\beta_2$, we conclude from \eqref{4.1.4} and \eqref{4.1.4a} that for any $\gamma > \gamma_1$ we can fix $\alpha >1$ satisfying \eqref{4.1.8} and \eqref{4.1.81} and obtain
	$$\frac{d}{dt} \left( \io (u+1)^\gamma + \io |\nabla v|^{2\alpha} \right)
	\le c_3 \io (u+1)^\gamma + c_6$$
	for all $t \in (0, \tm)$. Thus, Gronwall's inequality implies that \eqref{wzor} holds for any
	$\gamma \in (\gamma_1, \infty)$. This proves the claim as $u$ is also uniformly bounded in $L^1 (\Omega)$.
\qed

{\bf Proof of Theorem \ref{theo2}.} Due to (\ref{wzor}) and the classical regularity theory of parabolic equations
applied to the second equation of \eqref{0}, see \cite[Lemma 4.1]{horstmann_winkler} for example, one obtains an estimate of $\nabla v$ in $L^\infty (\Omega \times (0,T))$ for any finite $T \in (0, \tm]$. Hence, we are now in a position to apply
\cite[Lemma A.1]{taowin_jde} in order to gain an estimate of $u$ in $L^\infty
(\Omega \times (0,T))$ which shows the existence of a global solution. More precisely, keeping the notation of \cite[Lemma A.1]{taowin_jde}, we have $f:= -\psi(u)\nabla v$ and $g:=0$, while due to (\ref{balance})
we can choose $m=1-p$ and make sure that $\psi(u)$ grows at most polynomially with respect to $u$. Moreover, by
Lemma~\ref{lemat}, we have $u \in L^\infty ((0,T); L^{p_0} (\Omega))$ and $f \in L^\infty ((0,T); L^{q_1} (\Omega))$ for any $p_0 \in (1,\infty)$ and $q_1 \in (1,\infty)$. This freedom of choosing any $p_0<\infty$ as well as any $q_1<\infty$ enables us to make sure that all the assumptions of  \cite[Lemma A.1]{taowin_jde} are satisfied.

Furthermore, if we additionally assume that \eqref{G1} and \eqref{H1} are satisfied in the case $n \ge 3$ and \eqref{GH2}
holds in the case $n=2$, we apply
\cite[Theorem~5.1]{win_mmas} in order to deduce that $(u,v)$ blows up in infinite time for suitably chosen initial data.
This finishes the proof of Theorem~\ref{theo2}.
\qed

%\begin{rem}
%As it was recently noticed in \cite{ishida} Winkler's strategy of using \eqref{convex} to obtain energy estimates for chemotaxis systems can work also in non-convex domains. It seems that following the approach in \cite{ishida} one can remove the assumption of convexity of the domain in Lemma \ref{lemat} and consequently in Theorem~\ref{theo2}.  
%\end{rem}

\mysection{Applications to volume filling models}\label{section5}

The last section is devoted to the analysis of the so-called volume filling models introduced in \cite{hil:volume}.
These are the models of chemotactic movement taking into account the size of cells. Since the size (though being small)
is nonzero, a probability that a cell attains a position in a neighborhood of a point $(x,t)$ depends on the density of cells at $(x,t)$. One of the models proposed in \cite{hil:volume} consists of the system \eqref{0} with
\begin{equation}\label{erst}
\phi(u):=(1+u)^{-\gamma}+u\gamma(1+u)^{-\gamma-1} ,\;\;\psi(u):=u(1+u)^{-\gamma}, \qquad \gamma >0.
\end{equation}
In what follows we state a theorem on different kind of behaviors of solutions to \eqref{0} with nonlinear diffusion and sensitivity given by \eqref{erst} in dimensions $n\geq 3$ for different $\gamma$. Moreover, we shall show that in
a space dimension $2$ for any $\gamma > 0$ volume filling models admit critical mass, i.e. there exists such a value of initial mass, that if a solution starts with a mass smaller than it, it exists globally as a bounded solution, while if the initial mass exceeds the above mentioned critical value, solution becomes unbounded in infinite time (still existing for any time $t>0$). One-dimensional problem has been studied in \cite{zhang_zheng} and it was proved that bounded solutions exist globally in time.

As a consequence of Theorems~\ref{theo1} and \ref{theo2} we have the following
\begin{prop}\label{prop1}
Consider the problem \eqref{0}, with $\phi$ and $\psi$ given by \eqref{erst}, in a bounded domain with a smooth boundary. In dimensions $n\geq 3$ 
%for $\gamma<1-\frac{2}{n}$ there exist radially symmetric solutions blowing up at a finite time $T<\infty$
%(this is not proved so far, as we have $p=\gamma >0$ and $q = 1-\gamma \in (\frac{2}{n},1)$), while 
for $\gamma>2-\frac{2}{n}$ all the solutions to \eqref{0}, \eqref{erst} exist globally in time. Moreover in the case $\gamma>2-\frac{2}{n}$, there are radially symmetric solutions blowing up at infinity.
\end{prop}

Quite an interesting situation appears in dimension $n=2$ where the critical mass phenomenon is observed. Namely we have the following
\begin{theo}\label{theo3}
Consider the problem \eqref{0}, with $\phi$ and $\psi$ given by \eqref{erst}, in $\Omega\subset \R^2$ with a smooth boundary and
nonnegative initial data $u_0 \in C(\bar{\Omega})$ and  $v_0\in W^{1,\infty} (\Omega)$. 
\begin{enumerate}
  \item[(i)] If $\int_\Omega u_0(x) dx< 4\pi (1+\gamma)$, then the unique solution to \eqref{0}, \eqref{erst} is global in time and 
	  bounded.  
	\item[(ii)] Assume further that $u_0$ and $v_0$ are radially symmetric. Then in case of $\int_\Omega u_0(x) dx< 8\pi (1+\gamma)$ the 
	  unique solution to \eqref{0}, \eqref{erst} is global in time and bounded. Moreover, for $\int_\Omega u_0(x) dx >
		8\pi (1+\gamma)$ %, not equal to $8k(1+\gamma)\pi$, $k\in \N$, 
		there is a solution which becomes unbounded in finite or infinite time.
\end{enumerate}	
\end{theo}

%\begin{rem}
%In the case of radially symmetric solutions, with $\Omega=B(0,R)$, a bounded global solution exists for any initial data $u_0$ such that $\int_\Omega u_0(x) dx< 8\pi(1+\gamma)$, while in the condition guaranteeing infinite-time blowup, one needs to replace $4$ by $8$. During the proof we will make comments related to those issues.
%\end{rem}

\begin{rem}
According to Theorem~\ref{theo2} we know that if $\Omega$ is bounded with smooth boundary then for $\gamma>1$ each solution to \eqref{0}, \eqref{erst}, with no restriction on initial data, exists globally in time. In particular, in view of Theorem \ref{theo3}, we know that for $\gamma>1$ unbounded radially symmetric solutions blow up at infinity.
\end{rem}

Our proof uses the ideas which appeared in the semilinear case, see \cite{B, NSY} for global existence and \cite{Ho02} for unboundedness. However, at some points where it seems that the straightforward application of methods from \cite{B, NSY, Ho02} is not possible, we provide necessary modifications. Our proof consists of several steps, a first one, where one notices what is the value of a critical mass, is the following lemma.

\begin{lem}\label{lem4.0}
Consider the problem \eqref{0}, with $\phi$ and $\psi$ given by \eqref{erst}, in $\Omega\subset \R^2$ with a smooth boundary. For a nonnegative $u_0 \in C(\bar{\Omega})$ such that $\int_\Omega u_0(x) dx< 4\pi (1+\gamma)$ and any nonnegative $v_0\in W^{1,\infty}(\Omega)$
(or $\int_\Omega u_0(x) dx< 8\pi (1+\gamma)$ in the radially symmetric setting) the Liapunov functional ${\cal F}$ is bounded from below and there exists such a positive constant $c>0$ that
\begin{equation}\label{lem4.0.1}
\sup_{0\le t<\tm} \max\{\int_\Omega|\nabla v|^2 dx, \int_\Omega u\ln u dx\}< c.
\end{equation}
Moreover for any $t \le \tm$ there exists a positive constant $c_1$ such that
\begin{equation}\label{4.0.4}
\int_0^t\int_\Omega v_t^2dxdt + \int_0^t\int_\Omega \frac{|\phi(u)\nabla u-\psi(u)\nabla v|^2}{\psi(u)}dxdt \le c_1.
\end{equation}
\end{lem}

\proof  Let us first compute precisely the Liapunov functional ${\cal F}$ in the case of $\phi$ and $\psi$ given by (\ref{erst}). In the definition of $G$ we take $s_0:=1$ and see that
\begin{equation}\label{G.0}
{\cal F}(u, v)=\int_\Omega G(u)dx -\int_\Omega uv dx + \frac{1}{2} \int_\Omega |\nabla v|^2 + \frac{1}{2} \int_\Omega v^2,
\end{equation}
where
\begin{equation}\label{G.1}
G(u):=\gamma (1+u)\ln (1+u)+u\ln u-u\left(\gamma(1+\ln 2)+1\right)-\gamma\ln 2+ \gamma+1.
\end{equation}
Hence
\begin{equation}\label{laposz}
\frac{{\cal F}}{\gamma+1}\geq \int_\Omega u\ln udx-\frac{1}{\gamma+1}\int_\Omega uvdx+\frac{1}{2(\gamma +1)}\int_\Omega |\nabla v|^2 dx
+\frac{1}{2(\gamma+1)}\int_\Omega v^2dx+C(m,\gamma, |\Omega|).
\end{equation}
Next, we use the idea in \cite{B, NSY}. We denote by
\[
\tilde{v}:=m\left(\int_\Omega e^{\frac{v}{\gamma+1}}dx\right)^{-1}e^{\frac{v}{\gamma+1}},\;\; m:=\int_\Omega udx.
\]
To shorten the formulas, let us moreover denote $\tilde{\mu}:=\int_\Omega e^{\frac{v}{\gamma+1}}dx$.

By the convexity of $\tau\rightarrow -\ln \tau$ we have
\[
0=-\ln \left(\left\|u\right\|_{L^1(\Omega)}^{-1}\int_\Omega u\left(\frac{\tilde{v}}{u}\right)\right)\leq \left\|u\right\|_{L^1(\Omega)}^{-1}\int_\Omega u\ln\left(\frac{u}{\tilde{v}}\right),
\]
hence
\[
0\leq \int_\Omega u\ln\left(\tilde{\mu}m^{-1}ue^{-\frac{v}{\gamma+1}}\right)
=m\ln \tilde{\mu}-m\ln m+\int_\Omega u\ln udx-\frac{1}{\gamma+1}\int_\Omega uv.
\]
Using (\ref{laposz}) we arrive at
\begin{equation}\label{4.0.0}
0\leq m\ln \tilde{\mu}+\frac{{\cal F}}{\gamma+1}-\frac{1}{2(\gamma+1)}\int_\Omega |\nabla v|^2 dx-\frac{1}{2(\gamma+1)}\int_\Omega v^2dx+C(m,\gamma, |\Omega|).
\end{equation}
Next, we are in a position to apply the Trudinger-Moser inequality (see \cite{NSY} and the references therein) to the function $v$
\[
\ln \left(\int_\Omega e^{\frac{v}{\gamma+1}}\right)\leq \left(C+\frac{1}{(\gamma+1)^2\beta}\int_\Omega |\nabla v|^2\right),
\]
where $\beta\in (0,8\pi)$, however in the case of radially symmetric functions in a ball $\beta\in (0,16\pi)$.
Inserting the above inequality into \eqref{4.0.0} we arrive at
\[
0\leq \left(\frac{m}{\beta(1+\gamma)^2}-\frac{1}{2(\gamma+1)}\right)\int_\Omega |\nabla v|^2+\frac{{\cal F}}{\gamma+1}-\frac{1}{2(\gamma+1)}\int_\Omega v^2dx+C(m,\gamma, \int_\Omega vdx, |\Omega|).
\]
Hence for $m < \frac{(\gamma+1)\beta}{2}$,  we arrive at
the following bound for $\nabla v$
\begin{equation}\label{4.0.1}
\int_\Omega|\nabla v|^2dx\leq C\left[\frac{{\cal F}(u_0,v_0)}{\gamma+1}+C(m,\gamma, \int_\Omega vdx, |\Omega|)\right].
\end{equation}
and the lower bound for the Liapunov function ${\cal F}$
\begin{equation}\label{4.0.2}
{\cal F}\geq -C(m,\gamma,\int_\Omega vdx, |\Omega|).
\end{equation}
Observe that the condition $m < \frac{(\gamma+1)\beta}{2}$ means $m < 4\pi (\gamma +1)$ (or
$m < 8\pi(\gamma+1)$ in the radial setting).
Next, Young's inequality gives
\[
\int_\Omega uvdx\leq \frac{1}{2} \int_\Omega u\ln u-\frac{m}{2}\ln 2+\int_\Omega e^{2v-1},
\]
which in turn gives
\[
{\cal F}(u_0,v_0)\geq {\cal F}(u,v)\geq  \frac{1}{2} \int_\Omega u\ln udx-\int_\Omega e^{2v-1}+\frac{1}{2}\int_\Omega |\nabla v|^2dx+\frac{1}{2}\int_\Omega v^2dx +C,
\]
using Trudinger-Moser's inequality and (\ref{4.0.1}), we thus obtain
\begin{equation}\label{4.0.3}
{\cal F}(u_0,v_0) +C\geq  \frac{1}{2}\int_\Omega u\ln udx.
\end{equation}
Moreover, \eqref{4.0.2} together with \eqref{liapunov}, \eqref{D} yields (\ref{4.0.4})
for any $0<t\leq \tm$. 
\qed

Next we introduce a definition of a set of stationary radially symmetric solutions to \eqref{0} and two lemmas which are the key steps in proving the unboundedness part of Theorem \ref{theo3}.
\begin{defi}\label{defiS}
Let $\Omega=B(0,R)\subset \R^2$. We say that nonnegative radially symmetric functions $(\bar{u},\bar{v}):\Omega\rightarrow \R^2$ belong to the set ${\cal S}$ if  for some constant $d\in \R$  $\bar{u}=\Xi^{-1}(\bar{v}+d)$, $\io \bar{u}dx=\io \bar{v}dx=m$ and $\bar{v}$ satisfies the following boundary value problem
\begin{equation}\label{A.1}
\Delta \bar{v}- \bar{v} + \Xi^{-1}(\bar{v}+d)=0\;\;\mbox{in}\;\;\Omega,\;\;\frac{\partial \bar{v}}{\partial\nu}=0\;\;\mbox{on}\;\;\partial\Omega,
\end{equation}
for $\Xi(s)$ being a primitive of $\frac{\phi(s)}{\psi(s)}$.
\end{defi}

\begin{lem}\label{lem4.1.1}
Assume that $u$ solving \eqref{0} with radially symmetric initial conditions $(u_0, v_0)$ is bounded in $L^\infty(\Omega)$. Then there exists $(\bar{u},\bar{v})\in {\cal S}$ such that for some subsequence of time $t_k \rightarrow\infty$ $(u(t_k,\cdot),v(t_k,\cdot))\rightarrow (\bar{u},\bar{v})$ in $C(\bar{\Omega})\times C^1(\bar{\Omega})$.
\end{lem}

\proof We start by applying the regularity theory of parabolic equations to the lower equation in \eqref{0} to arrive at the uniform-in-time estimate of $\|v(t, \cdot)\|_{C^\alpha(\bar{\Omega})}$. Next, since we are in the case of a parabolic system with a triangular main part, we can apply classical theory of parabolic systems, see for instance \cite{amann}, to find a bound independent of time of the %$C^\alpha(\bar{\Omega})\times C^{1,\alpha}(\bar{\Omega})$ 
$C^{2,\alpha}(\bar{\Omega})\times C^{2,\alpha}(\bar{\Omega})$ norm of the couple $(u,v)$. Thus, we are allowed to apply the Arzel\`{a}-Ascoli theorem and extract a subsequence of times $t_k$ along which $(u(t_k,\cdot),v(t_k,\cdot))\rightarrow (\bar{u},\bar{v})$ to some $(\bar{u},\bar{v})$ in $C(\bar{\Omega})\times C^1(\bar{\Omega})$. It is enough if we show that $(\bar{u},\bar{v})$ belongs to ${\cal S}$ defined in Definition \ref{defiS}. To this end we make use of the LaSalle principle and the entropy production terms ${\cal D}$, see (\ref{D}) for the definition, in the following way. The functions $(\bar{u},\bar{v})$ are such that the functional ${\cal F}$ when evaluated on the trajectory of solutions to \eqref{0} starting from $(\bar{u},\bar{v})$ is constant. This means that ${\cal D}$ is $0$ when evaluated at this trajectory. Hence $(\bar{u},\bar{v})\in {\cal S}$.
\qed

The following lemma shows that one can choose radially symmetric initial data such that the value of the Liapunov functional ${\cal F}(u_0,v_0)$ is arbitrarily small.

\begin{lem}\label{lem4.1.2}
Assume $m>8\pi (1+\gamma)$ and consider the functional ${\cal F}$ given by \eqref{G.0}, \eqref{G.1} with $\Omega=B(0,R) 
\subset \R^2$.
There exists a sequence of nonnegative radially symmetric functions $(u_k, v_k)$ satisfying $\io u_kdx=m$ such that
\begin{equation}\label{4.1.3a}
{\cal F}(u_k, v_k)\rightarrow -\infty.
\end{equation}
\end{lem}

\proof In the proof we will use the ideas from \cite{struwe}. First we bound the Liapunov functional ${\cal F}$ from above making use of \eqref{G.1} and the inequality $x\ln x\leq (x+1)\ln (x+1)$ for $x>0$
\begin{eqnarray}\label{4.6.1}
\frac{{\cal F}(u,v)}{(\gamma+1)}-c(m, \gamma, |\Omega|)&\leq &\nn \\ 
\io (u+1)\ln(u+1)&-&\frac{1}{\gamma+1}\io (u+1)v+ \frac{1}{\gamma+1}\io v+\frac{1}{2(\gamma+1)}\io |\nabla v|^2+\frac{1}{2(\gamma+1)}\io v^2.\nn \\
\end{eqnarray}
Obviously, unboundedness from below of the right-hand side of \eqref{4.6.1} means that ${\cal F}$ is unbounded. 
Next choose $u$ in the form 
\begin{equation}\label{4.6.15}
u+1:=\frac{(m+|\Omega|)e^{\frac{v}{1+\gamma}}}{\io e^{\frac{v}{1+\gamma}}dx}\;.
\end{equation}
In view of 
\begin{equation}\label{4.6.16}
\frac{1}{\gamma+1}\io (u+1)vdx=\frac{m+|\Omega|}{(\gamma+1)\io e^{\frac{v}{1+\gamma}}dx}  \io v e^{\frac{v}{1+\gamma}}dx
\end{equation}
and 
\begin{eqnarray}\label{4.6.17}
&&\io (u+1)\ln (u+1)dx\nn\\
&=&\frac{m+|\Omega|}{\io e^{\frac{v}{1+\gamma}}dx}\left(\ln (m+|\Omega|)\io e^{\frac{v}{1+\gamma}}dx-\io e^{\frac{v}{1+\gamma}}\left(\ln\io e^{\frac{v}{1+\gamma}}dx\right)dx+ \io e^{\frac{v}{1+\gamma}}\frac{v}{1+\gamma}dx\right),\nn \\
\end{eqnarray}
\eqref{4.6.1} yields
\begin{eqnarray}\label{4.6.18}
&&{\cal F}\left(\frac{(m+|\Omega|)e^{\frac{v}{1+\gamma}}}{\io e^{\frac{v}{1+\gamma}}dx}\;
,\; v\right)-c(m,\gamma, |\Omega|)\nn \\
&\leq & -(m+|\Omega|)(1+\gamma)\ln \io e^{\frac{v}{1+\gamma}}dx + \io v+\frac{1}{2}\io |\nabla v|^2+\frac{1}{2}\io v^2.
\end{eqnarray}
Next, we notice that by Jensen's inequality
\begin{equation}\label{4.6.19}
-|\Omega|\ln \io e^{\frac{v}{1+\gamma}}dx\leq -|\Omega|\ln |\Omega|-\frac{1}{\gamma+1}\io vdx.
\end{equation}
Hence in the light of \eqref{4.6.18} and \eqref{4.6.19} proving the lemma is reduced to finding a sequence of nonnegative radially symmetric functions $v_k$ such that for $m > 8\pi (1+ \gamma)$
\begin{equation}\label{4.6.2}
\frac{1}{2}\io |\nabla v_k|^2dx+\frac{1}{2}\io v_k^2dx-m(\gamma+1)\ln\io e^{\frac{v_k}{1+\gamma}}dx\rightarrow -\infty 
\end{equation}
when $k\rightarrow \infty$ (in particular notice that a sequence $u_k$ associated to $v_k$ by \eqref{4.6.15} satisfies $\io u_kdx=m$).

Substituting $z:=\frac{v}{1+\gamma}$ in \eqref{4.6.2} we notice that finding a sequence of radially symmetric functions $z_k>0$ such that the functional
\begin{equation}\label{4.6.3}
F(z):=-\frac{m}{1+\gamma}\ln \io e^zdx+\frac{1}{2}\io |\nabla z|^2dx+\frac{1}{2}\io z^2dx
\end{equation} 
goes to $-\infty$ when evaluated on $z_k$ is enough.

To this end we notice first that if we find a sequence $z_k$, not necessarily positive, such that $F(z_k)\rightarrow -\infty$, still $z_k^+$ is a sequence of nonnegative functions such that $F(z_k^+)\rightarrow -\infty$. As $z_k$ we take an example from \cite{struwe}. Namely, 
\[
z_k:=\ln \frac{(1/k)^2}{((1/k)^2+\pi|x|^2)^2}- \frac{1}{|\Omega|}\io \ln \frac{(1/k)^2}{((1/k)^2+\pi|x|^2)^2}dx.
\]
In \cite{struwe} it was shown that a functional very similar to $F$ goes down to $-\infty$ when evaluated on $z_k$, for reader's convenience we provide an argument that it is also the case for $F$ (a similar argument can be found in \cite{wang_wei}).

First
\[
\io |\nabla z_k|^2dx=32\pi\ln k+O(1), \io \frac{(1/k)^2}{((1/k)^2+\pi|x|^2)^2}dx=O(1),
\]
where $O(1)\leq C$ when $k\rightarrow \infty$. Next
\begin{equation}\label{zetka}
\frac{1}{|\Omega|} \io \ln \left(\frac{(1/k)^2}{((1/k)^2+\pi|x|^2)^2}\right)dx=2\ln 1/k+O(1). 
\end{equation}
Moreover, by \eqref{zetka} we have
\[
\io z_k^2dx=\io\left(\ln\frac{(1/k)^2}{((1/k)^2+\pi|x|^2)^2}-\ln(1/k)^2+O(1)\right)^2dx,
\]
\[
\io\left(\ln \frac{1}{((1/k)^2+\pi|x|^2)^2}\right)^2dx=O(1).
\]
Summing up all the above calculations, we arrive at the following estimate
\[
F(z_k)=\left(16\pi-\frac{2m}{1+\gamma}\right)\ln k+O(1),
\]
and the lemma follows.
\qed

As a further step towards the proof of Theorem~\ref{theo3} we prove the following result, see \cite{zhang_zheng}.
\begin{prop}\label{prop5}
  Consider the problem \eqref{0}, with $\phi$ and $\psi$ given by \eqref{erst}, in $\Omega\subset \R^2$ with a smooth boundary.
	In case of $\int_\Omega u_0(x) dx< 4\pi (1+\gamma)$ (or $\int_\Omega u_0(x) dx< 8\pi (1+\gamma)$ in the radially symmetric setting)
	and $\gamma>1$, for any $t_0 \in (0,\tm)$ there is a constant $c_2(t_0)>0$ such that
	\begin{equation}\label{eq5}
	  \left( \io v_t^2 (x,t) \, dx \right)^{\frac{1}{2}} \le c_2(t_0), \qquad t \in (t_0, \tm). 
	\end{equation}
\end{prop}

\proof Observe that by parabolic regularity the solution $(u,v)$ to \eqref{0} is smooth in $\Omega \times (0,\tm)$. Differentiating the 
  second equation of \eqref{0} with respect to time, multiplying by $v_t$ and integrating with respect to $x$, we arrive at
	\begin{eqnarray*}
	  \frac{1}{2} \frac{d}{dt} \io v_t^2 &=& - \io |\nabla v_t|^2 - \io v_t^2 + \io u_t v_t \nn \\
		&=& - \io |\nabla v_t|^2 - \io v_t^2 - \io \left( \phi(u) \nabla u - \psi (u) \nabla v \right) \nabla v_t \nn \\
		&\le& - \frac{1}{2} \io |\nabla v_t|^2 - \io v_t^2 + \frac{1}{2} \io \left| \phi(u) \nabla u - \psi (u) \nabla v \right|^2, 
	\end{eqnarray*}
  where in the first line we used that $|\nabla v_t| =0$ on $\partial \Omega$ holds, which is a consequence of the Neumann boundary
	condition. Next observe that $\gamma >1$ implies $\psi \le 1$ so that
	$$\frac{1}{2} \io \left| \phi(u) \nabla u - \psi (u) \nabla v \right|^2 \le 
	  \frac{1}{2} \io \frac{\left| \phi(u) \nabla u - \psi (u) \nabla v \right|^2}{\psi(u)}.$$
	Hence, we obtain
	$$\frac{1}{2} \frac{d}{dt} \io v_t^2 + \frac{1}{2} \io |\nabla v_t|^2 + \io v_t^2	
	  \le \frac{1}{2} \io \frac{\left| \phi(u) \nabla u - \psi (u) \nabla v \right|^2}{\psi(u)}, \qquad t \in (0, \tm).$$
	Finally integrating the latter with respect to $t$ and	using \eqref{4.0.4} we conclude that
	$$ \io v_t^2 (x,t) \, dx \le 2c_1 + \io v_t^2(x,t_0) \, dx <\infty, \qquad t \in (t_0,\tm),$$
	is satisfied for any $t_0 \in (0,\tm)$.
\qed

The final preparation for the proof of Theorem~\ref{theo3} is the following lemma.
\begin{lem}\label{lem5.1}
   Consider the problem \eqref{0}, with $\phi$ and $\psi$ given by \eqref{erst}, in $\Omega\subset \R^2$ with a smooth boundary. 
	 If $\int_\Omega u_0(x) dx< 4\pi (1+\gamma)$ (or $\int_\Omega u_0(x) dx< 8\pi (1+\gamma)$ in the radially symmetric setting) 
	 then there
	 exists $C>0$ such that 
   \begin{equation}\label{5.1.6}
     \| u(t)\|_{L^{\gamma +2}(\Omega)} \le C \qquad\mbox{ for all } t \in (0, \tm).
   \end{equation}	
\end{lem}

\proof For any $\gamma>0$, multiplying the first equation of \eqref{0} by $(1+u)^{\gamma+1}$ and using the second equation of \eqref{0}, 
we obtain
\begin{eqnarray}\label{5.1.1}
  & & \hspace*{-20mm} \frac{1}{\gamma +2} \frac{d}{dt} \io (1+u)^{\gamma+2} \nn \\
	&=& - (\gamma+1) \io \phi(u) (1+u)^\gamma |\nabla u|^2 + (\gamma+1) \io \psi(u) (1+u)^\gamma \nabla u \nabla v \nn \\
	&\le& - (\gamma+1) \io |\nabla u|^2 + (\gamma+1) \io u \nabla u \nabla v \nn \\
	&=& - (\gamma+1) \io |\nabla u|^2 + \frac{\gamma+1}{2} \io \nabla(u^2) \nabla v \nn \\
	&=& - (\gamma+1) \io |\nabla u|^2 - \frac{\gamma+1}{2} \io u^2 (v_t +v-u) \nn \\
	&\le& - (\gamma+1) \io |\nabla u|^2 + \frac{\gamma+1}{2} \left( \io u^4 \right)^{\frac{1}{2}} \left( \io v_t^2 \right)^{\frac{1}{2}} 
	+ \frac{\gamma+1}{2} \io u^3.
\end{eqnarray}
In view of \eqref{lem4.0.1}, we are in a position to apply \cite[(22)]{debye} and for any $\eps >0$ arrive at
\begin{equation}\label{5.1.2}
  \| u\|_{L^3 (\Omega)}^3 \le \eps \|\nabla u\|_{L^2(\Omega)}^2 + C(\eps, m, \| u \ln u\|_{L^1 (\Omega)}).
\end{equation}
Now we treat cases $\gamma>1$ and $\gamma\le 1$ separately. In the case $\gamma\le 1$ the proof follows the lines of 
\cite{NSY}. The case  $\gamma>1$ requires another idea which is based on Proposition~\ref{prop5}.

\noindent
{\bf Case $\gamma>1$.} In view of the inequalities of Gagliardo-Nirenberg and Young, for any $\eps >0$ we have
\begin{equation}\label{5.2.1}
  \| u\|_{L^4 (\Omega)}^2 \le C \| u \|_{H^1 (\Omega)}^{\frac{3}{2}} 
	\|u \|_{L^1(\Omega)}^{\frac{1}{2}} \le \eps \| \nabla u\|_{L^2(\Omega)}^2 + C(\eps,m).
\end{equation} 
Using once more the Gagliardo-Nirenberg inequality and defining $\alpha := \frac{2(\gamma+2)}{\gamma+1}$, we obtain
$$\|1+u \|_{L^{\gamma+2} (\Omega)}^{\alpha} \le C \left( \| \nabla u\|_{L^2 (\Omega)}^{\frac{\alpha(\gamma +1)}{\gamma+2}} 
  \|1+u \|_{L^1 (\Omega)}^{\frac{\alpha}{\gamma+2}}
  + \|1+u \|_{L^1 (\Omega)}^{\alpha} \right) \le c_3 \left( 1+ \| \nabla u\|_{L^2 (\Omega)}^2 \right),$$
where $c_3 = c_3 (m)$.		
This implies
\begin{equation}\label{5.2.2}
  \| \nabla u\|_{L^2 (\Omega)}^2 \ge \frac{1}{c_3} \|1+u \|_{L^{\gamma+2} (\Omega)}^{\alpha} -1.
\end{equation}	
Defining $t_0 := \frac{\tm}{2}$, $\eps := \frac{1}{1+c_2(t_0)}$, combining \eqref{5.1.1}-\eqref{5.2.2} and applying 
Proposition~\ref{prop5} (which requires $\gamma >1$), we deduce that
\begin{eqnarray*}
  \frac{1}{\gamma +2} \frac{d}{dt} \io (1+u)^{\gamma+2} 
	&\le& - \frac{\gamma+1}{2} \io |\nabla u|^2 + C(t_0, m, \| u \ln u\|_{L^1 (\Omega)}) \nn \\
	&\le& - \frac{\gamma+1}{2 c_3} \|1+u \|_{L^{\gamma+2} (\Omega)}^{\alpha} + C(t_0, m, \| u \ln u\|_{L^1 (\Omega)})
\end{eqnarray*}
for $t \in (t_0, \tm)$. In view of $\alpha>0$ this implies
$$\io (1+u)^{\gamma+2} (x,t) \, dx \le \max \left\{ \io (1+u)^{\gamma+2} (x,t_0) \, dx, \left( \frac{2c_3}{\gamma +1} 
   C(t_0, m, \| u \ln u\|_{L^1 (\Omega)}) \right)^{\frac{\gamma+2}{\alpha}} \right\}$$
for $t \in (t_0, \tm)$.	According to the local existence result, we conclude that \eqref{5.1.6} holds.

\noindent
{\bf Case $\gamma\le 1$.} 
Using the Gagliardo-Nirenberg and Young inequalities and defining $\theta := \frac{2-\gamma}{4} \in (0,\frac{1}{2})$, we estimate
the second term on the right-hand side of \eqref{5.1.1} for any $\eps >0$:
\begin{eqnarray}\label{5.1.3}
  \|v_t\|_{L^2 (\Omega)} \| u\|_{L^4 (\Omega)}^2 &\le& C \|v_t\|_{L^2 (\Omega)} \| u \|_{H^1 (\Omega)}^{2\theta} 
	\|u \|_{L^{\gamma+2}(\Omega)}^{2(1-\theta)} \nn \\
	&\le& \eps \| \nabla u\|_{L^2(\Omega)}^2 + \eps \| u\|_{L^2(\Omega)}^2 + C_\eps \|v_t\|_{L^2 (\Omega)}^{\frac{1}{1-\theta}} 
	\|u \|_{L^{\gamma+2}(\Omega)}^2 \nn \\
	&\le& 2\eps \| \nabla u\|_{L^2(\Omega)}^2 + C_\eps \| u\|_{L^1(\Omega)}^2 + C_\eps \left( \|v_t\|_{L^2 (\Omega)}^2 + 1 
	\right) \|u \|_{L^{\gamma+2}(\Omega)}^2 \nn \\
	&\le& 2\eps \| \nabla u\|_{L^2(\Omega)}^2 + \eps \left( \|v_t\|_{L^2 (\Omega)}^2 +1\right) \|u \|_{L^{\gamma+2}(\Omega)}^{\gamma +2} 
	\nn \\ 
	& & + C(\eps,m) \left( \|v_t\|_{L^2 (\Omega)}^2 +1\right).
\end{eqnarray} 
Inserting \eqref{5.1.2} and \eqref{5.1.3} into \eqref{5.1.1} and using $u \ge 0$, for any $\eps \in (0, \frac{1}{6})$  we arrive at
\begin{eqnarray}\label{5.1.4} 
  \frac{1}{\gamma +2} \frac{d}{dt} \io (1+u)^{\gamma+2}
  &\le& - \frac{\gamma+1}{2} \io |\nabla u|^2 + \frac{\gamma+1}{2} \left( \|v_t\|_{L^2 (\Omega)}^2 + \eps\right) \|1+ u \|_{L^{\gamma+2}(\Omega)}^{\gamma +2} \nn \\
	& & + C(\eps,m, \| u \ln u\|_{L^1 (\Omega)}) \left( \|v_t\|_{L^2 (\Omega)}^2 +1\right). 
\end{eqnarray}
In view of the Gagliardo-Nirenberg inequality and $\gamma \le 1$, there exists $C_0 = C_0 (m)>0$ such that
$$\|1+u \|_{L^{\gamma+2} (\Omega)}^{\gamma +2} \le C \left( \| \nabla u\|_{L^2 (\Omega)}^{\gamma +1} \|1+u \|_{L^1 (\Omega)}
  + \|1+u \|_{L^1 (\Omega)}^{\gamma +2} \right) \le C_0 \left( 1+ \| \nabla u\|_{L^2 (\Omega)}^2 \right),$$
which implies
$$
  \| \nabla u\|_{L^2 (\Omega)}^2 \ge \frac{1}{C_0} \|1+u \|_{L^{\gamma+2} (\Omega)}^{\gamma +2} -1.
$$	
Inserting this into \eqref{5.1.4}, we obtain
\begin{eqnarray}\label{5.1.5}
   \frac{1}{\gamma +2} \frac{d}{dt} \io (1+u)^{\gamma+2}
  &\le& \frac{\gamma+1}{2} \left( \|v_t\|_{L^2 (\Omega)}^2 + \eps - \frac{1}{C_0} \right) \io (1+u)^{\gamma +2} \nn \\
	& & + C(\eps,m, \| u \ln u\|_{L^1 (\Omega)}) \left( \|v_t\|_{L^2 (\Omega)}^2 +1\right). 
\end{eqnarray}
Choosing now $\eps \in (0, \min \{\frac{1}{6}, \frac{1}{2C_0} \})$ and defining $y(t) := 	\io (1+u)^{\gamma+2}$, by \eqref{5.1.5} we get
$$y^\prime (t) \le \frac{(\gamma+2)(\gamma+1)}{2} \left( \|v_t (t)\|_{L^2 (\Omega)}^2 - \frac{1}{2C_0} \right) y(t) + C_1
  \left( \|v_t (t)\|_{L^2 (\Omega)}^2 +1\right), \qquad t \in (0,\tm).$$
Setting $a(t) := \int_0^t \frac{(\gamma+2)(\gamma+1)}{2} \left( \|v_t (s)\|_{L^2 (\Omega)}^2 - \frac{1}{2C_0} \right) \, ds$,
this yields
\begin{eqnarray*}
  y(t) &\le& y(0) e^{a(t)} + C_1 \int_0^t \left( \|v_t (s)\|_{L^2 (\Omega)}^2 +1\right) e^{a(t)-a(s)} \, ds.
\end{eqnarray*} 
In view of \eqref{4.0.4}, there are positive constants $C_2$ and $C_3$ such that
$$a(t) -a(s) \le C_2 - C_3 (t-s) \qquad \mbox{for all } 0 \le s \le t < \tm$$
is fulfilled. Hence, we obtain
\begin{eqnarray*}
  y(t) \le y(0) e^{C_2 - C_3 t} + C_1 e^{C_2} \int_0^{\tm} \io v_t^2 (x,s) dx ds + C_1 e^{C_2} \frac{1-e^{-C_3 t}}{C_3}, \qquad
	t \in (0, \tm).	
\end{eqnarray*}	
This proves \eqref{5.1.6}. 
\qed

\noindent
\vspace{0.3cm}
{\bf Proof of Theorem \ref{theo3}.}

The proof splits into two parts. In the first of them we show the unboundedness above the critical mass. The second one is devoted to proving that global solutions are bounded for initial mass less than $4\pi (1+\gamma)$, $\gamma > 0$
(or for the initial mass less than $8\pi (1+\gamma)$ in the case of radially symmetric initial data).

First part follows the strategy in \cite{Ho02}. First notice that due to uniqueness of local solutions and rotational invariance of the operators in both equations of \eqref{0}, if we start from radially symmetric initial data, this property of solution is preserved for any $t>0$. In order to show that solutions starting from radially symmetric initial data with mass larger than $8\pi (1+\gamma)$ are unbounded, we assume the contrary. Next, using Lemma \ref{lem4.1.1} and continuity of the Liapunov functional we infer that ${\cal F}(u_0(x),v_0(x))\geq \inf_{(\bar{u},\bar{v})\in {\cal S}} {\cal F}(\bar{u},\bar{v})$.  However, Lemma \ref{Appendix} tells us that there exists a constant $\bar{c}\in \R$ such that
\[
\inf_{(\bar{u},\bar{v})\in {\cal S}} {\cal F}(\bar{u},\bar{v})> \bar{c},
\]
which is in contradiction with Lemma \ref{lem4.1.2}, since according to this lemma we can pick up such initial data (radially symmetric) that 
\[
{\cal F}(u_0(x),v_0(x))< \bar{c}.
\]

Next we proceed with a proof of the claim concerning solutions with initial mass less than $4\pi (1+\gamma)$
(or less than $8\pi (1+\gamma)$ in the radially symmetric setting). As $\gamma >0$ and $n=2$,
using Lemma~\ref{lem5.1} and the classical regularity theory of parabolic equations
applied to the second equation of \eqref{0}, see \cite[Lemma 4.1]{horstmann_winkler} for example, we obtain a constant $C>0$ such that 
\begin{equation}\label{5.3.1}
  \| \nabla v \|_{L^\infty (\Omega \times (0,\tm))} \le C.
\end{equation}	
This enables us to provide for any $\alpha \in  (\gamma +2,3\gamma +4]$ the uniform estimate
\begin{equation}\label{5.3.2}
  \| u(t) \|_{L^\alpha (\Omega)} \le C(\alpha), \qquad\mbox{for all } t \in (0,\tm).
\end{equation}
To this end we fix $\alpha \in  (\gamma +2,3\gamma+4]$, multiply the first equation of \eqref{0} by $(1+u)^{\alpha-1}$ and use 
\eqref{5.3.1} to arrive at
\begin{eqnarray}\label{5.3.2a}
  \frac{1}{\alpha} \frac{d}{dt} \io (1+u)^\alpha 
	&\le& - (\alpha-1) \io (1+u)^{\alpha-\gamma-2} |\nabla u|^2 + (\alpha -1) \io (1+u)^{\alpha -\gamma -1} |\nabla u \nabla v| \nn \\
	&\le& - \frac{\alpha-1}{2} \io (1+u)^{\alpha-\gamma-2} |\nabla u|^2 + C_1 \io (1+u)^{\alpha -\gamma} \nn \\
	&=& - \frac{2(\alpha -1)}{(\alpha -\gamma)^2} \io \left| \nabla (1+u)^{\frac{\alpha - \gamma}{2}} \right|^2 + C_1 \io (1+u)^{\alpha -
	\gamma}.
\end{eqnarray}
Next, our aim is to show that 
\begin{equation}\label{5.3.3}
  \io (1+u)^{\alpha - \gamma} \le \frac{\alpha -1}{C_1(\alpha - \gamma)^2}\io \left| \nabla (1+u)^{\frac{\alpha - \gamma}{2}} \right|^2 
	+ C_2 \left( \io (1+u)^{\gamma +2} \right)^{\frac{\alpha -\gamma}{\gamma +2}}. 
\end{equation}
To this end we define $f := (1+u)^{\frac{\alpha - \gamma}{2}}$ and $\beta := \frac{2(\gamma+2)}{\alpha -\gamma}$. 
First we notice that when $\beta \ge 2$, \eqref{5.3.3} holds due to H\"{o}lder's inequality. The case when 
$\beta \in [1,2)$ requires some more effort. The 
Gagliardo-Nirenberg inequality yields
$$\|f \|_{L^2(\Omega)}^2 \le C_3 \left( \|\nabla f \|_{L^2(\Omega)}^{2 \theta} \|f \|_{L^\beta (\Omega)}^{2(1-\theta)} 
  + \|f \|_{L^\beta (\Omega)}^2 \right) 
	\le \frac{\alpha -1}{C_1(\alpha - \gamma)^2} \|\nabla f \|_{L^2(\Omega)}^2 + C_2 \|f \|_{L^\beta (\Omega)}^2$$  
with $\theta := 1 - \frac{\beta}{2} \in (0,1)$. In view of the choices of $f$ and $\beta$, this implies \eqref{5.3.3}. Using once more the Gagliardo-Nirenberg inequality, we obtain 
$$\|f \|_{L^\frac{2\alpha}{\alpha -\gamma} (\Omega)} \le C \left( \|\nabla f \|_{L^2(\Omega)}^{\mu} \|f \|_{L^\beta (\Omega)}^{1-\mu} 
  + \|f \|_{L^\beta (\Omega)} \right)$$
with $\mu := \frac{\alpha - \gamma -2}{\alpha} \in (0,1)$ due to $\alpha > \gamma +2$. In view of Lemma~\ref{lem5.1}, this implies
\begin{equation}\label{5.3.4}
  \io \left| \nabla (1+u)^{\frac{\alpha - \gamma}{2}} \right|^2 \ge C \left( \io (1+u)^\alpha	\right)^{\frac{\alpha-\gamma}{\alpha \mu}}
		-1.
\end{equation}	
By \eqref{5.3.2a}-\eqref{5.3.4}, there exist positive constants $C_4$ and $C_5$ such that
$$\frac{d}{dt} \io (1+u)^\alpha \le C_4 - C_5 \left( \io (1+u)^\alpha	\right)^{\frac{\alpha-\gamma}{\alpha - \gamma -2}},
 \qquad t \in (0, \tm).$$
This proves \eqref{5.3.2} for any $\alpha \in (\gamma +2, 3\gamma +4]$.

Now we are in a position to apply \cite[Lemma A.1]{taowin_jde} in order to gain an estimate of $u$ in $L^\infty
(\Omega \times (0,\tm))$ which shows that $(u,v)$ is global in time and bounded. In fact, keeping the notation of \cite[Lemma A.1]{taowin_jde} with $f:= -\psi(u)\nabla v$ and $g:=0$, in view of \eqref{erst} we can fix $m=1-\gamma$, $q_1 := 3\gamma+ 4$ and 
$p_0 := 3\gamma +4$. Then, by \eqref{5.3.1} and
\eqref{5.3.2}, we have $u \in L^\infty ((0,T); L^{p_0} (\Omega))$ and $f \in L^\infty ((0,T); L^{q_1} (\Omega))$. Moreover, the above
choices of $p_0$ and $q_1$ imply $q_1 >n+2$ and $p_0 > \max \{ \gamma, 2\gamma -1, 4, 3\gamma + \frac{2}{3}\}$ so that all the assumptions of 
\cite[Lemma A.1]{taowin_jde} are satisfied.
This completes the proof of the theorem.
\qed

\mysection{Appendix}\label{section6}

Appendix is devoted to proving Lemma \ref{Appendix} which states the bound from below of the Liapunov functional over the set of stationary solutions to the 2d volume filling Keller-Segel system with a power-type probability jump function.  Our approach is based on \cite{HW}, however we have to deal with a more general system of stationary solutions due to the fact that the function $G$ appearing in the Liapunov functional is more general in our case. We restrict ourselves to the radially symmetric case, but our proof seems to be more straightforward and self-contained. We only refer the reader to the famous Br\'ezis-Merle inequality, \cite{Brezis_Merle}, and the rest of the argument is presented in detail.

\begin{lem}\label{Appendix}
For $m\neq 8\pi(1+\gamma)$ %q$, $q\in \N$, 
the values of ${\cal F}$ over the set ${\cal S}$ are bounded from below.
\end{lem}
\proof
Let us recall that ${\cal F}$ is given by 
\[
\io G(u)dx-\io uvdx+\frac{1}{2}\io |\nabla v|^2dx+\frac{1}{2}\io v^2dx
\] 
with $G$ being defined in \eqref{G.1}. Next, we recall that ${\cal S}$ is defined in Definition \ref{defiS} as
radially symmetric functions satisfying \eqref{A.1} in a ball $B(0,R)$.

Our proof of Lemma \ref{Appendix} will follow by contradiction. If $v$ belongs to ${\cal S}$, then $u=\Xi^{-1}(v+d)$ as was noticed in Definition~\ref{defiS}. Let us first assume that there exists a sequence of 
functions $v_k\in C^1(\bar{B(0,R)})\cap {\cal S}$, $k\in \N $, such that $\io u_kdx=m$ and real numbers $d_k$, $k\in \N $, such that 
\begin{equation}\label{App1}
\left\|\nabla v_k\right\|_{L^2(B(0,R))}\rightarrow \infty,
\end{equation}
\begin{equation}\label{App2}
d_k\rightarrow -\infty
\end{equation}
and
\begin{equation}\label{App3}
\max_{x\in \bar{B(0,R)}}v_k(x)\rightarrow \infty.
\end{equation}
Let us notice that if at least one of \eqref{App1}, \eqref{App2} or \eqref{App3} is not true than the claim of Lemma \ref{Appendix} follows. Indeed, it is clear that if there exists a constant $C>0$ such that $\max_{x\in \bar{B(0,R)}}v_k(x)\leq C$ then ${\cal F}$ is bounded over ${\cal S}$. Next let us show that denying \eqref{App1} also leads to the claim of Lemma \ref{Appendix}. 
%Indeed, by \eqref{4.0.0} we estimate ${\cal F}$ from below by 
%\[
%-(\gamma+1)m\ln \io e^{\frac{v_k}{1+\gamma}}dx+\frac{1}{2}\int_\Omega |\nabla v_k|^2 dx+\frac{1}{2}\int_\Omega v_k^2dx+C(m,\gamma, |\Omega|).
%\]
%By the Trudinger-Moser inequality, since \eqref{App1} is violated and we have a bound of $\left\|\nabla v_k\right\|_{L^2(B(0,R))}$, we see that there exists $C>0$ such that $\ln \io e^{\frac{v_k}{1+\gamma}}dx\leq C$, so 
%${\cal F}$ is bounded.
Indeed, by \eqref{A.1} and the Gagliardo-Nirenberg inequality, we obtain
$$\int_\Omega u_k v_k dx = \int_\Omega | \nabla v_k |^2 dx + \int_\Omega v_k^2 dx \le 2 \int_\Omega | \nabla v_k |^2 dx
+ c(\Omega) \| v_k \|_{L^1(\Omega)}^2.$$
Since $\int_\Omega v_k = m$ and the fact that \eqref{App1} is violated and we have a bound of $\left\|\nabla v_k\right\|_{L^2(B(0,R))}$, we see that there exists $C>0$ such that $\io u_k v_k dx\leq C$, so that ${\cal F}$ is bounded from below. 

Finally, we prove also that when \eqref{App2} is violated then ${\cal F}$ is bounded.  To this end we notice that up to a constant $u\Xi(u)\approx G(u)$, so that 
\[
{\cal F}(u_k,v_k)\geq \io u_k\left(\Xi(u_k)-v_k+const\right)+\frac{1}{2}\left\|v_k\right\|_{H^1(\Omega)}^2+c(m,\gamma, \Omega),
\] 
and since $(u_k,v_k)\in {\cal S}$ we see that by Definition \ref{defiS} 
\begin{equation}\label{5.0.63}
u_k(\Xi(u_k)-v_k)=d_ku_k, 
\end{equation}
what gives
\[
{\cal F}(u_k,v_k)\geq d_km +\frac{1}{2}\left\|v_k\right\|_{H^1(\Omega)}^2+c(m,\gamma, \Omega)\geq C
\]
if \eqref{App2} does not hold. In the next step of the proof we notice that since $\io u_kdx=m$, by Prokhorov's theorem
we extract a subsequence (still denoted by $u_k$) such that 
\begin{equation}\label{A.2}
u_k\stackrel{*}\rightharpoonup \mu
\end{equation} 
in measures. We will show that \eqref{App1}-\eqref{App3} yields
\begin{equation}\label{5.0.1}
\io u_k dx\rightarrow 8\pi(1+\gamma), %q \;\;\mbox{for some}\;\; q\in \N,
\end{equation}
this would contradict the assumption of Lemma \ref{Appendix}. 

Let us define the set of blowup points
\begin{equation}\label{BS}
{\cal BS}:=\{x\in \bar{\Omega}: \exists x_k\stackrel{k\rightarrow \infty}\longrightarrow x\;\mbox{such that}\;\;  v_k(x_k)\stackrel{k\rightarrow \infty}\longrightarrow \infty \}.
\end{equation}
Next, for $\delta>0$ let us define a $\delta$-regular point in the following way (see \cite{Brezis_Merle}, \cite{wang_wei}), $x_0$ is $\delta$-regular if there exists 
compactly supported smooth function $0\leq \zeta\le 1$ with $\zeta=1$ in some neighborhood of $x_0$ such that
\begin{equation}\label{reg}
\io \zeta d\mu<\frac{3\pi(1+\gamma)}{1+2\delta}.
\end{equation} 
By $\Sigma(\delta)$ we denote the set of points which are not $\delta$-regular.

The next step of the proof consisits of the following two propositions, their proofs we postpone till the end of the section. 
\begin{prop}\label{dreg}
If $x_0$ is $\delta$-regular then there exists $r>0$ such that the sequence $v_k$ solving \eqref{A.1} is bounded in the ball $B(x_0,r)$ uniformly in $k\in \N$. Moreover, ${\cal BS}=\Sigma(\delta)$ and ${\cal BS}$ contains only finitely many elements. 
If $v_k$ are radially symmetric then ${\cal BS}=\{0\}$. 
\end{prop}
In view of Proposition \ref{dreg} we see that ${\cal BS}=\{0\}$. Thus we define 
\[
\sigma^k(r):=\int_{B(0,r)}u_kdx.
\]
We have the following fact.
\begin{prop}\label{mu-char}
Let $u_k, v_k\in {\cal S}$, in particular $v_k, d_k$ solve \eqref{A.1}, then the limit measure (see \eqref{A.2}) $\mu=c\delta_0$, where $c>0$ and 
\begin{equation}\label{mu_wzor}
\lim_{k\rightarrow \infty}\io u_kdx=\lim_{r\rightarrow 0}\lim_{k\rightarrow \infty}\sigma^k(r).
\end{equation}
\end{prop}
In the sequel we shall need the following fact, an argument validating it can be found in \cite[p.233]{wang_wei} or \cite[p.168]{HW}.
\begin{prop}\label{5.0.3}
Let $\Omega\subset \R^2$ and $f_k\in C^1(\bar{\Omega})$ be a radially symmetric function such that $\left\|f_k\right\|_{L^1(\Omega)}\leq C$ for some positive $C$ and all $k \in \N$.
If $v_k$ satisfies 
\begin{equation}\label{5.0.35}
-\Delta v_k+ v_k =f_k\;\;\mbox{in}\;\;B(0,R),\;\;\frac{\partial v_k}{\partial\nu}=0\;\;\mbox{on}\;\;\partial B(0,R),
\end{equation} 
then for small $0<r<R$ and $0\neq|x|\leq r$ we have
\begin{equation}\label{5.0.4}
v_k(x)=-\frac{1}{2\pi}\ln |x|\int_{B(0,r)}f_k(y)dy+C(x,r,f_k)
\end{equation}
and
\begin{equation}\label{5.0.5}
\nabla v_k(x)=-\frac{1}{2\pi}\frac{x}{|x|^2}\int_{B(0,r)}f_k(y)dy+C(x,r,f_k),
\end{equation}
where $C(x,r,f_k)$ is such that $C(x,r,f_k)\stackrel{r\rightarrow 0}\longrightarrow 0$ uniformly in $k\in \N$ and $|x|\leq r$.
\end{prop} 
Next we proceed to show that,  
\begin{equation}\label{5.0.2}
\lim_{r\rightarrow 0}\lim_{k\rightarrow \infty}\sigma^k(r)=8(1+\gamma)\pi
\end{equation}
what yields \eqref{5.0.1} and we arrive at a contradiction.
Since $v_k\in {\cal S}$ is radially symmetric, it satisfies \eqref{A.1}. 
In radial coordinates \eqref{A.1} reads
\[
(r(v_k)_r)_r= -r\Xi^{-1}(v_k+d_k)+ rv_k, \;(v_k)_r(0)=(v_k)_r(R)=0.
\]
We multiply it by $(v_k)_r$, integrate over $B(0,r)$ and arrive at the following Pohozaev identity 
\begin{equation}\label{5.0.6}
\pi \left((v_k)_rr\right)^2-\pi r^2v_k^2+2\pi r^2F(v_k)=2\int_{B(0,r)}F(v_k)dx-\int_{B(0,r)}v_k^2dx,
\end{equation}
where $F(s):=\int_0^s\Xi^{-1}(\sigma+d_k)d\sigma$. We notice that in view of Definition \ref{defiS}, \eqref{erst} and
$\Xi(s) = \ln s + \gamma \ln (s+1)$, we have
\begin{equation}\label{5.0.62a}
e^{\frac{v_k}{1+\gamma}}e^{\frac{d_k}{1+\gamma}}-1\leq u_k = \Xi^{-1}(v_k+d_k)\leq e^{\frac{v_k}{1+\gamma}}e^{\frac{d_k}{1+\gamma}}
\end{equation}
and so
\begin{equation}\label{5.0.62}
(1+\gamma)\left( e^{\frac{v_k}{1+\gamma}}-1 \right)e^{\frac{d_k}{1+\gamma}}-v_k\leq F(v_k)\leq (1+\gamma) \left( e^{\frac{v_k}{1+\gamma}}-1 \right)e^{\frac{d_k}{1+\gamma}}\leq (1+\gamma)e^{\frac{d_k}{1+\gamma}}e^{\frac{v_k}{1+\gamma}}.
\end{equation}
Using Proposition \ref{5.0.3} we can estimate some terms appearing in \eqref{5.0.6}, namely  
\begin{equation}\label{5.0.8}
r^2v_k^2(r)\leq Cr 
\end{equation}
by \eqref{5.0.4}, while \eqref{5.0.62} leads to
\begin{equation}\label{5.0.85}
r^2F(v_k(r))\leq r^2(1+\gamma)Ce^{\frac{d_k}{1+\gamma}}e^{-\ln (r^{\frac{\sigma^k(r)}{2\pi(1+\gamma)}})}= r^2(1+\gamma)C\;
\frac{e^{\frac{d_k}{1+\gamma}}}{r^{\frac{\sigma^k(r)}{2\pi(1+\gamma)}}},
\end{equation}
where $\sigma^k(r) =\int_{B(0,r)}u_k dx \in [0,m]$, so for any $r>0$
\begin{equation}\label{5.0.87}
r^2F(v_k(r))\rightarrow 0\;\;\mbox{when}\;\;k\rightarrow \infty
\end{equation}
in view of \eqref{App2}.

Since $u_k$ is estimated in $L^1$, by elliptic regularity results we have for $1\leq q<2$
\begin{equation}\label{5.0.61}
\left\|\nabla v_k\right\|_{L^q(B(0,r)}\leq C,
\end{equation}
hence 
\begin{equation}\label{5.0.7}
\int_{B(0,r)}v_k^2dx\leq Cr^{1/2} \left\|v_k\right\|_{W^{1,\frac{3}{2}}(\Omega)}^2\leq Cr^{1/2}.
\end{equation}
Next we observe that by \eqref{5.0.62a} and \eqref{5.0.62}
\[
(1+\gamma) \left( u_k (r) - e^{\frac{d_k}{1+\gamma}} \right) -v_k(r)\leq F(v_k(r))\leq 
(1+\gamma)(u_k (r) +1).
\]
Hence, using also \eqref{5.0.7} we obtain
\begin{equation}\label{5.0.10}
(1+\gamma)\sigma^k(r)-O(r^2) e^{\frac{d_k}{1+\gamma}}-O(r^{1/4})\leq \int_{B(0,r)}F(v_k)dx \leq (1+\gamma)\sigma^k(r)+ O(r^2).
\end{equation}
Moreover, by \eqref{5.0.5} we have
\begin{equation}\label{5.0.9}
\pi(v_k)_r^2(r)r^2=\left(\frac{\sigma^k(r)}{2\pi}\right)^2\left(\pi +O(r)\right).
\end{equation}
We plug \eqref{5.0.8}, \eqref{5.0.87}, \eqref{5.0.7}-\eqref{5.0.9} in \eqref{5.0.6} and run the two scale argument. First we let 
$k\rightarrow\infty$, next we take the limit when $r\rightarrow 0$. We have 
\[
2(1+\gamma)\lim_{r\rightarrow 0}\lim_{k\rightarrow \infty}\sigma^k(r)=\frac{\pi}{4\pi^2}\left(\lim_{r\rightarrow 0}\lim_{k\rightarrow \infty}\sigma^k(r)\right)^2.
\]
Thanks to the above identity and in view of \eqref{mu_wzor} we arrive at
\[
m=\lim_{r\rightarrow 0}\lim_{k\rightarrow \infty}\sigma^k(r)=8\pi(1+\gamma),
\]
a contradiction.
\qed
Now, in order to complete the argument, we prove Propositions \ref{dreg} and \ref{mu-char}.

\vspace{0.3cm}
{\bf Proof of Proposition \ref{dreg}}. Take $x_0$ to be $\delta$-regular point. First we notice that for small $0<\rho$ and $1<q<2$ 
\begin{equation}\label{5.5.1}
\int_{B(x_0,\rho)}v_kdx\leq C\left\|v_k\right\|_{L^q(B(0,R))}\rho^{\frac{2(q-1)}{q}}.
\end{equation} 
Hence, in view of \eqref{5.0.61}, for arbitrary
$\varepsilon>0$ we can choose $\rho$ small enough to make sure that $\int_{B(x_0,\rho)}v_kdx<\varepsilon$. 
Next, we split $v_k$ into two parts $v_k:=v_{1k}+v_{2k}$, where $v_{ik}, i=1,2$ satisfy
\begin{eqnarray}\label{5.5.3}
-\Delta v_{1k}&=&u_k-v_k\;\;\mbox{in}\;\; B(x_0,\rho),\\
v_{{1k}_{|\partial B(x_0,\rho)}}&=&0\nn    
\end{eqnarray}
and
\begin{eqnarray}\label{5.5.4}
-\Delta v_{2k}=0\;\;\mbox{in}\;\; B(x_0,\rho),\\
v_{{2k}_{|\partial B(x_0,\rho)}}=v_{{k}_{|\partial B(x_0,\rho)}},\nn
\end{eqnarray}
respectively. We notice that by the maximum principle $v_{2k}>0$. Since $x_0$ is $\delta$-regular, by \eqref{reg}, for $k$ large enough we have
\begin{equation}\label{5.5.45}
\int_{B(x_0,\rho)}u_k\leq \frac{3\pi(1+\gamma)}{1+\frac{3}{2}\delta}\;,
\end{equation}
Next, we recall the Br\'ezis-Merle inequality (see \cite[Theorem 1]{Brezis_Merle}), then for any $b\in (0,4\pi)$
\begin{equation}\label{5.5.2}
\int_{B(x_0,\rho)} e^{\frac{(4\pi-b)v_{1k}}{\left\|u_k-v_k\right\|_{L^1(B(x_0,\rho))}}}dx\leq C.
\end{equation}
Matching \eqref{5.5.45} with \eqref{5.5.1} yields
\begin{equation}\label{5.5.5}
\int_{B(x_0,\rho)}|u_k-v_k| dx\leq \frac{3\pi(1+\gamma)}{1+\delta}\;.
\end{equation}
Our next observation is that $\int_{B(x_0,\rho)}v_{2k}dx\leq\int_{B(x_0,\rho)}\left(v_k+|v_{1k}|\right)dx\leq C$, 
and by the Harnack inequality we see
\begin{equation}\label{5.5.6}
\left\|v_{2k}\right\|_{L^\infty(B(x_0,\rho/2))}\leq C\left\|v_{2k}\right\|_{L^1(B(x_0,\rho/2))}\leq C.
\end{equation}
On the other hand, by \eqref{5.5.5} and \eqref{5.5.2} with $b=\pi$ we arrive at
\begin{equation}\label{5.5.7}
\int_{B(x_0,\rho)}e^{\frac{v_{1k}(1+\delta)}{1+\gamma}}dx\leq \int_{B(x_0,\rho)}e^{\frac{3\pi v_{1k}}{\left\|u_k-v_k\right\|_{L^1(B(x_0,\rho))}}}dx\leq C.
\end{equation}
As a consequence of \eqref{5.5.7} we obtain $\left\|u_k\right\|_{L^{1+\delta}(B(x_0,\rho))}\leq C$, and consequently
by \eqref{5.5.6} 
\[
\left\|u_k-v_k\right\|_{L^{1+\delta}(B(x_0,\rho/2))}\leq C
\] 
due to \eqref{5.0.61}, this in turn implies by standard elliptic regularity (notice that $v_{1k}$ satisfies \eqref{5.5.3}) that 
\[
\left\|v_{1k}\right\|_{L^\infty(B(x_0,\rho/4))}\leq C.
\]
and the boundedness claim of Proposition \ref{dreg} holds.  

Now we pass to the second claim of Proposition \ref{dreg}. First we handle the inclusion ${\cal BS}\subset \Sigma(\delta)$. Indeed, if $x_0\in {\cal BS} \setminus \Sigma(\delta)$, then it is $\delta$-regular
and in the light of just proven boundedness part of Proposition \ref{dreg} we arrive at a contradiction.

To see the opposite inclusion we notice that if $x_0\in \Sigma(\delta)$ then for any $\rho>0$
\[
\left\|v_k\right\|_{L^\infty(B(x_0,\rho))}=\infty.
\]
Indeed, otherwise there exists $\rho_0>0$ and a subsequence $v_k$ such that all $k\in \N$
\[
\left\|v_k\right\|_{L^\infty(B(x_0,\rho_0))}\leq C,
\] 
consequently
\[
e^{\frac{d_k}{1+\gamma}}e^{\frac{v_k}{1+\gamma}}\leq C e^{\frac{d_k}{1+\gamma}}
\]
in $B(x_0, \rho_0)$. But since $d_k\rightarrow -\infty$ we arrive at a contradiction, $x_0\notin \Sigma(\delta)$.

Next we notice that $x_0\in \Sigma(\delta)$ if and only if
$\mu(\{x_0\})\geq \frac{3\pi(1+\gamma)}{1+\delta}$. Hence the cardinality of $\Sigma(\delta)$ is finite, 
precisely speaking less or equal than $\frac{m(1+\delta)}{3\pi(1+\gamma)}$. In particular, in the case of radially
symmetric solutions $\Sigma(\delta)=\{0\}$.
\qed

\vspace{0.3cm}
{\bf Proof of Proposition \ref{mu-char}}. It is sufficient to show that $\mu=c\delta_0$ for some $c>0$. The measure $\mu$ is the weak star limit of $u_kdx$, while $0\leq u_k\leq e^{d_k}e^{v_k}$. Outside any ball $B(0,r)$ functions $v_k$ are bounded uniformly in $k$ by Proposition~\ref{dreg}. Hence $0\leq u_k\leq Ce^{d_k}$, since $d_k\rightarrow -\infty$ we see that $u_k$ tends to $0$ uniformly in $k$ outside a ball of radius $r>0$ for any ball. This means that $u_k\stackrel{*}\rightharpoonup c\delta_0$.  
\qed

%\vspace{0.3cm}
%{\bf Proof of Proposition \ref{5.0.3}}. Functions $v_k$ can be represented with the use of Green's functions
%as 
%\[
%v_k(x)=\int_{B(0,R)} G(x,y)f_k(y)dy,
%\]
%where $G(x,y)=-\frac{1}{2\pi} \ln |x-y|+H(x,y)$, $H$ is a harmonic correction.  To see \eqref{5.0.4} it is enough to estimate the asymptotic behavior for $r\approx 0$ of 
%\[
%-\frac{1}{2\pi} \int_{B(0,R)} \ln |x-y|f_k(y)dy,
%\]
%to this end we notice that
%\begin{eqnarray}\label{5.1.01}
%\int_{B(0,R)} \ln |x-y|f_k(y)dy=\int_{B(0,r)}\ln \left(|x-y|-|x|\right)f_k(y)dy \nn \\
%+\ln |x|\int_{B(0,r)}f_k(y)dy+\int_{B(0,R)\setminus B(0,r)} \ln |x-y|f_k(y)dy.
%\end{eqnarray}
%In view of the uniform in $k$ bound on $f_k$ in $L^1$, and owing Young's inequality, last term on the right-hand side 
%of \eqref{5.1.01} is estimated by the constant independent of $r$ and $k$. Next, for small $\eta>0$ we estimate
%\begin{eqnarray}\label{5.1.02}
%\int_{B(0,r)}\ln \left(|x-y|-|x|\right)f_k(y)dy=\int_{\{|x-y|\geq \eta\}}\ln \left(|x-y|-|x|\right)f_k(y)dy\nn \\
%+\int_{\{|x-y|< \eta\}}\ln \left(|x-y|-|x|\right)f_k(y)dy
%\end{eqnarray}
%The proof of \eqref{5.0.5} goes along the same reasoning.

%\qed

\noindent
{\bf Acknowledgement.} Both authors are grateful to Michael Winkler from Paderborn for his suggestion to use his technique to prove the
global existence part in this paper and for his help concerning the borderline case $p=0$ of the finite-time blowup. T.C. was partially supported by the National Centre of Science (NCN) under grant 2013/09/D/ST1/03687. T.C. is grateful to Sasha Mikhaylov from the Steklov Institute for helpful discussion. This work
was initiated during a visit of C.~Stinner at the Instytut Matematyczny PAN in Warsaw. He is grateful for the invitation, support and
hospitality.

\end{document}